\documentclass[11pt]{amsart} 
\usepackage[a4paper]{geometry}

\usepackage[all]{xy}
\usepackage{verbatim}
\usepackage{amsmath,amsthm,amssymb,amscd}
\usepackage{tikz}
\usepackage{epsfig,amssymb,latexsym,amscd,amsmath,amsthm,verbatim}
\usepackage{bbm}
\usepackage[mathcal]{euscript}
\usepackage{appendix}

\input{grcawol.sty}
\input grcalcx.sty

\theoremstyle{plain}

\newtheorem{teo}{Teorema}[section]
\newtheorem{thm}[teo]{Theorem}
\newtheorem{cor}[teo]{Corollary}
\newtheorem{lemma}[teo]{Lemma}
\newtheorem{prop}[teo]{Proposition}

\theoremstyle{remark}

\newtheorem{rema}[teo]{Remark}

\theoremstyle{definition}

\newtheorem{defi}[teo]{Definition}

\newtheorem{ex}[teo]{Example}

\newcommand{\BB}{{\mathcal{B}}}
\newcommand{\CC}{{\mathcal{C}}}
\newcommand{\DD}{{\mathcal{D}}}

\newcommand{\ACA}{{{}_A\mathcal{C}_A}}

\newcommand{\HYD}{{}^H_H\mathcal {YD}}
\newcommand{\kGYD}{{}^{\Bbbk G}_{\Bbbk G}\mathcal {YD}}
\newcommand{\kZYD}{{}^{\Bbbk \Z}_{\Bbbk \Z}\mathcal {YD}}

\newcommand{\AkGYDA}{{\,_A\!\left(\kGYD\right)\!_A}}

\newcommand{\id}{\ensuremath{\mathrm{id}}}

\newcommand{\Hom}{\ensuremath{\mathrm{Hom}}}
\newcommand{\Homk}{\ensuremath{\mathrm{Hom}}_\Bbbk}

\newcommand{\HomAA}{\Hom_{AA}}

\newcommand{\homAA}{\hom_{AA}}

\newcommand{\Vect}{\ensuremath{\mathrm{Vect}}}
\newcommand{\Mod}{\ensuremath{\mathrm{Mod}}}
\newcommand{\AModA}{{_A\Mod_A}}
\newcommand{\Ch}{\ensuremath{\mathrm{Ch}}}
\newcommand{\Coch}{\ensuremath{\mathrm{Coch}}}
\newcommand{\hCh}{{\overline{\Ch}}}
\newcommand{\hCoch}{{\overline{\Coch}}}
\newcommand{\khCh}{{\underline{\Ch}}}
\newcommand{\khCoch}{{\underline{\Coch}}}

\newcommand{\Ext}{\ensuremath{\mathrm{Ext}}}

\newcommand{\ant}{s}

\newcommand{\N}{{\mathbb N}}
\newcommand{\Z}{{\mathbb Z}}

\renewcommand{\H}{\ensuremath{\mathrm{H}}}

\newcommand{\op}{\ensuremath{\mathrm{op}}}

\newcommand{\flop}{{\op}}

\newcommand{\totimes}{{\;\Tilde{\otimes}\;}}
\newcommand{\totimesA}{{\;\Tilde{\otimes}_A\;}}
\newcommand{\todot}{{\;\Tilde{\odot}\;}}
\newcommand{\tdiamond}{{\;\Tilde{\diamond}\;}}
\newcommand{\tstar}{{\;\Tilde{\star}\;}}
\newcommand{\AW}{\ensuremath{\mathrm{AW}}}
\newcommand{\AWt}{\overline{\ensuremath{\mathrm{AW}}}}

\newcommand{\cc}{\left(}
\newcommand{\dd}{\right)}

\newcommand{\opo}{\otimes\cdots\otimes}
\newcommand{\inv}{{{}^{-1}}}

\title{Graded braided commutativity in Hochschild cohomology}

\author{Javier C\'oppola}
\address{Instituto de Matem\'atica y Estad\'\i stica ``Rafael Laguardia'', Facultad de Ingenier\'\i a, Universidad de la Rep\'ublica, Uruguay}
\email{jcoppola@fing.edu.uy}

\author{Andrea Solotar}
\address{IMAS-CONICET and Departamento de Matem\'atica, Facultad de Ciencias Exactas y Naturales, Universidad de Buenos Aires, Argentina }
\email{asolotar@dm.uba.ar}
\date{}

\thanks{\footnotesize This work has been supported by the projects UBACYT 20020170100613BA, PIP-
CONICET 11220200101855CO and Mathamsud-AREPTHEO The first mentioned author was supported by doctoral scholarships from CAP-UdelaR. The second
mentioned author is a research member of CONICET (Argentina), Senior Associate at ICTP
and visiting Professor at Guangdong Technion-Israel Institute of Technology.}

\begin{document}
\maketitle

\begin{abstract}
We prove the graded braided commutativity of the Hochschild cohomology of $A$ with trivial coefficients, where $A$ is a braided Hopf algebra in the category of Yetter-Drinfeld modules over the group algebra of an abelian group, under some finiteness conditions on a projective resolution of $A$ as $A$-bimodule. This is a generalization of a result by Mastnak, Pevtsova, Schauenburg and Witherspoon to a context which includes Nichols algebras such as the Jordan and the super Jordan plane. We prove this result by constructing a coduoid-up-to-homotopy structure on the aforementioned projective resolution in the duoidal category of chain complexes of $A$-bimodules. We also prove that the Hochschild complex of a braided bialgebra $A$ in an arbitrary braided monoidal category is a cocommutative comonoid up to homotopy with the deconcatenation product which induces the cup product in Hochschild cohomology. \end{abstract}



\section{Introduction}
Given a field $\Bbbk$ and an associative algebra $A$, the Hochschild cohomology $\H^\bullet(A,A)$ has a very rich structure, which has been widely studied. In particular it is a nonnegatively graded algebra with the cup product which is also graded commutative; this fact allows using methods of commutative algebra. The graded commutativity has been first proved by Gerstenhaber \cite{gerstenhaber1964}, this proof gave rise to the Gerstenhaber bracket. Another proof of the graded commutativity was obtained via the Eckmann-Hilton argument \cite{marianoeckmannhilton}. In fact the cup product is also defined for $\H^\bullet(A,R)$, where $R$ is an $A$-bimodule with an associative product $R\otimes_A R\to R$, but graded commutativity is in general no longer true in this context. If $A$ is an augmented algebra, one can set $R=\Bbbk$, where $\Bbbk$ is an $A$-bimodule via the augmentation. It has been proven that for a Hopf algebra $A$, $\H^\bullet(A,\Bbbk)$ is also graded commutative \cite{farinatisolotar, marianoeckmannhilton,  taieffer2004}. The proof in \cite{farinatisolotar} uses the fact that $\H^\bullet(A,\Bbbk)$ is a subalgebra of $\H^\bullet(A,A)$. The question of what happens when $A$ is a Hopf algebra in a different category has been raised in \cite{mpsw}. In that article, the authors work with Hopf algebras in braided categories. They prove that in this situation, if certain internal hom objects exist, then $\H^\bullet(A,\Bbbk)$ is a graded braided commutative algebra. 
An example of this situation is considering Nichols algebras; these are Hopf algebras in the braided category of Yetter-Drinfeld modules over a Hopf algebra $H$. They have been thoroughly  studied in the case $H=\Bbbk G$ for an abelian or finite group $G$ - see for example \cite{andruskiewitsch2002pointed, aapw, stefanvay, estanislaofk2020}. The additional hypothesis about internal homs translates here to the condition of either $A$ or $H$ being finite dimensional.

Looking at the computations of the graded algebra structure of $\H^\bullet(A,\Bbbk)$ for $A$ the Jordan plane \cite{lopessolotar} or the super Jordan plane \cite{recasolotar}, we realized that even if the Hochschild cohomology algebra is not graded commutative, it is still graded braided commutative, even though the Jordan and the super Jordan plane are both infinite dimensional and the Hopf algebra is $\Bbbk\Z$, which is also infinite dimensional. The computations are detailed in the Appendix of this paper.

Our main aim is to prove that graded braided commutativity still holds in a more general context. For this, we show how to make use of arbitrary resolutions. The key step is moving from braided monoidal categories to duoidal categories.

The contents of this work are the following:

In Section 2 we recall some preliminary definitions and results, together with some examples that are relevant in the sequel: Hochschild cohomology in Subsection 2.1, monoidal and braided monoidal categories in Subsections 2.2 and 2.3 respectively, Nichols algebras in Subsection 2.4.

The subject of Section 3 is graded braided commutativity. We review some results and definitions from \cite{mpsw} and finally prove a theorem - see Theorem \ref{thmMPSWhomotopy} - which slightly generalizes Theorem 3.12 and Corollary 3.13 of \cite{mpsw}. An intermediate step for this result is a theorem of braided graded cocommutativity up to homotopy for the complex analog to the bar resolution - see Theorem \ref{thmdec} - which makes no use of internal homs.

Section 4 contains our main results about graded braided commutativity, Theorem \ref{thmmain} and Corollary \ref{cormain}. The cases of the Jordan and of the super Jordan plane are particular instances of these last results.

In this paper we will denote by $\Bbbk$ an arbitrary field. It will be clear from the context whether the symbol $\otimes$ represents the usual tensor product of $\Bbbk$-vector spaces, or stands for the product in an abstract monoidal category.

We thank Marcelo Aguiar and Estanislao Herscovich for fruitful discussions on this topic.

\section{Preliminaries}
\subsection{Hochschild cohomology}\label{sectionHochschild}

In this section we will recall some well-known facts about Hochschild cohomology that will be useful in this article, see for example \cite{witherspoonbook}.

\begin{defi}
Let $A$ be a $\Bbbk$-algebra and $M$ be an $A$-bimodule. The \emph{Hochschild cohomology of $A$ with coefficients in $M$} is defined as
\[
\H^\bullet(A,M) = \Ext^\bullet_{AA}(A,M),
\]
where $\Ext^\bullet_{AA}(A,-)$ denotes the left derived functors of $\Hom_{A\otimes A^\op}(A,-)$, that will be denoted $\HomAA(A,-)$.
\end{defi}

\begin{defi}
Consider $A$ as a bimodule over itself with actions given by its multiplication $\mu$. The \emph{bar resolution of $A$} is the following free resolution of $A$ as $A$-bimodule:
\[
\xymatrix{ \cdots \ar[r]^-{b'_3} & B_3(A) \ar[r]^-{b'_2}& B_2(A) \ar[r]^-{b'_1} & B_1(A) \ar[r]^-{b'_0} & B_0(A) \ar[r] & 0,
}
\]

 where $B_n(A) = A\otimes A^{\otimes n} \otimes A$ for all $n\in\N_0$ and the differential $b'_n: B_{n+1}(A)\to B_n(A)$ is given by the formula
 \[
 b'_n(a_0\otimes\cdots\otimes a_{n+2}) = \sum_{i=0}^{n+1} (-1)^i\ a_0\otimes\cdots\otimes a_i a_{i+1} \otimes\cdots\otimes a_{n+2}.
 \]
 The quasi-isomorphism from this resolution to $A$ is induced by $\mu:B_0(A)=A\otimes A\to A$.
\end{defi}

Thus, when computing Hochschild cohomology via this resolution one gets $\H^\bullet(A,M) = H(\Hom_{AA}(B_\bullet(A),M), (b'_\bullet)^*)$. The naturality of the adjunction $\Hom_{AA}(A\otimes V\otimes A,M) \simeq \Homk(V,M)$ implies that the complexes $(\Hom_{AA}(B_\bullet(A),M), (b'_\bullet)^*)$ and $(C^\bullet(A),d^\bullet)$ are isomorphic, where $C^n(A)= \Homk(A^{\otimes n},M)$ and $d^n:C^n(A)\to C^{n+1}(A)$ is:
\[
\begin{array}{rcl}
d^n f(a_1\otimes\cdots\otimes a_{n+1}) & = & a_1\cdot f(a_2\otimes\cdots\otimes a_{n+1})\\
&& +\ \sum_{i=1}^n (-1)^i\ f(a_1\otimes\cdots\otimes a_i a_{i+1}\otimes\cdots\otimes a_{n+1})\\
 && +\ (-1)^{n+1} f(a_1 \otimes\cdots\otimes a_n)\cdot a_{n+1}.
\end{array}
\]

\begin{rema} \label{remamultbimod}
Let $R$ be a $\Bbbk$-algebra, and $\rho: A\to R$ an algebra homomorphism. There is an $A$-bimodule structure on $R$ defining $a\cdot r\cdot b = \rho(a)r\rho(b),\forall a,b\in A,\ r\in R$. Moreover, with this bimodule structure the multiplication $\mu:R\otimes R\to R$ factors through $R\otimes_A R$.
\end{rema} 
\medskip
Let $P_\bullet \to A$ be a projective resolution of $A$ as $A$-bimodule. The canonical identification $A\overset{\simeq}\to A\otimes_A A$ lifts to a chain morphism $\omega:P_\bullet\to (P\otimes_A P)_\bullet$, unique up to a unique chain homotopy. This map allows the definition of a convolution product in $\Hom_{AA}(P_\bullet, R)$, denoted by $\cup$, which passes to cohomology. This is the \emph{cup product} of the Hochschild cohomology of $A$ with coefficients in $R$, denoted by $\cup: \H^\bullet(A,R) \otimes \H^\bullet(A,R) \to \H^\bullet(A,R)$.

For the bar resolution, one can explicitly define a coproduct as follows:
\[
\omega(1\otimes a_1\opo a_n\otimes 1) = \sum_{p+q=n} \cc 1\otimes a_1\opo a_p\otimes 1\dd\otimes_A \cc1\otimes a_{p+1}\opo a_{p+q}\otimes 1\dd,
\]
which gives the following well-known expression for the cup product:
\[
(\varphi\cup\psi)(1\otimes a_1\opo a_{p+q}\otimes 1) = \varphi(1\otimes a_1\opo a_p\otimes 1) \psi(1\otimes a_{p+1}\opo a_{p+q}\otimes 1),
\]
for all $\varphi\in\Hom_{AA}(B_p(A),R),\psi\in\Hom_{AA}(B_q(A),R),a_1,\ldots, a_{p+q}\in A$.

Translating this to the complex $(C^\bullet(A),d^\bullet)$ gives the following formula:
\[
(\alpha\cup\beta)(a_1\opo a_{p+q}) = \alpha(a_1\opo a_p)\beta(a_{p+1}\opo a_{p+q}),\ \forall a_1,\ldots,a_{p+q}\in A,
\]
$\alpha\in C^p(A),\beta\in C^q(A),a_1,\ldots a_{p+q}\in A$.

\begin{defi}
An \emph{augmented $\Bbbk$-algebra} is a $\Bbbk$-algebra $A$ equipped with an algebra homomorphism $\varepsilon:A\to\Bbbk$. We will omit $\varepsilon$ when it is clear from the context.
\end{defi}

Notice that if $A$ is an augmented $\Bbbk$-algebra then $\Bbbk$ becomes an $A$-bimodule. We will focus our attention in the algebra $\H^\bullet(A,\Bbbk)$.

Every Hopf algebra $A$ is an augmented algebra via its counit. Its cohomology with trivial coefficients is isomorphic to  $\Ext^\bullet_A(\Bbbk, \Bbbk)$, which is the so called \emph{Hopf algebra cohomology} \cite[Definition 9.3.5]{witherspoonbook}. Note that the category of $A$-bimodules is monoidal with the usual tensor product whenever $A$ is a Hopf algebra. After reviewing some concepts on braided monoidal categories, we will see a generalization of this fact in the last section.

\subsection{Monoidal categories}

In this section we will recall some facts about monoidal categories. For a more complete study in this topic we refer to \cite{maccategories} and \cite{aguiar}.

\begin{defi} A \emph{monoidal category} is a 6-tuple $(\CC, \otimes, I, a, \ell, r)$, where 

\begin{itemize}
\item $\CC$ is a category,
\item $\otimes: \CC \times \CC \to \CC$ is a bifunctor, called its \emph{product},
\item $I$ is an object in $\CC$, called its \emph{unit},
\item $a$ is a natural isomorphism, called \emph{associator}, given by a family of morphisms $a_{X,Y,Z}: (X \otimes Y) \otimes Z \to X \otimes (Y \otimes Z)$ for each $X,Y,Z\in\CC$, such that the following pentagonal diagram commutes for all $X,Y,Z,T\in\CC$: 
\[
\xymatrix{ & (X\otimes Y) \otimes (Z\otimes T)
\ar[dr]^{{a_{X,Y,Z\otimes T}}} & \\
 ((X\otimes Y) \otimes Z) \otimes T \ar[ur]^{{a_{X \otimes Y,Z,T}}} \ar[dd]_{a_{X,Y,Z}\otimes\id_T} &
& X\otimes (Y\otimes (Z\otimes T)) \\ \\
(X\otimes (Y \otimes Z)) \otimes T \ar[rr]_{a_{X,Y \otimes Z,T}} &
&X\otimes ((Y\otimes Z)\otimes T), \ar[uu]_{\id_X\otimes a_{Y,Z,T}}
}
\]
\item $\ell$ and $r$ are natural isomorphisms, called \emph{unitors}, given by families of morphisms $\ell_X: I\otimes X \to X$ and $r_X: X\otimes I \to X$  such that the following triangular diagram commutes for all $X,Y\in\CC$:
\[
\xymatrix{
(X\otimes I)\otimes Y \ar[rr]^{a_{X,I,Y}} \ar[dr]_{r_X \otimes \id_Y}
&& X\otimes (I\otimes Y) \ar[dl]^{\id_X \otimes \ell_Y} \\
& X\otimes Y. &
}
\]
\end{itemize}
A monoidal category is called \emph{strict} if $a, \ell, r$ are identities.

When the rest of the structure is clear from its context, we may refer to a monoidal category just by the name of the underlying category $\CC$. We will also say that it is a \emph{monoidal structure} on $\CC$.
\end{defi}
MacLane's Coherence Theorem \cite[Chapter VII, Section 2]{maccategories} states that these two relations imply that every composition of these morphisms from one arbitrary product of objects in $\CC$ to another one yields the same morphism.

\begin{rema}\label{remaflop}
If $(\CC, \otimes, I, a, \ell, r)$ is a monoidal category, then its \emph{opposite category} $\CC^\op$ admits a monoidal structure given by $(\CC^\op, \otimes, I, a^{-1}, \ell^{-1}, r^{-1})$.
\end{rema}

\begin{defi}\label{defitranspose}
Let $(\CC, \otimes, I, a, \ell, r)$ be a monoidal category. Its \emph{transpose} is the monoidal category $(\CC, \totimes, I, a^{-1}, r, \ell)$, where $X\totimes Y = Y\otimes X$ for all $X,Y\in\CC$.
\end{defi}

Without loss of generality, from now on we will consider strict monoidal categories, since every monoidal category is equivalent to a strict one via functors that preserve the monoidal structure \cite[Chapter XI, Section 3]{maccategories}, in a sense that will be precisely stated later. This is also the way monoidal categories are usually dealt with.

In this work, all monoidal categories will be $\Bbbk$-linear, being the product bilinear with respect to morphisms and distributive over direct sums. In order to have (co)chain complexes and their (co)homology inside these categories, we will need an abelian structure as well.

\begin{defi}
Let $(\CC, \otimes, I)$ be a monoidal category. A \emph{monoid} in $\CC$ is a triple $(A, \mu, \eta)$, where
\begin{itemize}
    \item $A$ is an object in $\CC$,
    \item $\mu:A\otimes A\to A$ is a morphism in $\CC$, called \emph{multiplication} or \emph{product}, such that the following diagram commutes:
    
    \[\xymatrix{
    A\otimes A\otimes A \ar[rr]^{\mu\otimes\id_A} \ar[d]_{\id_A\otimes\mu}& & A\otimes A \ar[d]^\mu \\
    A\otimes A \ar[rr]_\mu & & A,
    }\]
    
    This property is called \emph{associativity} of the product.
    
    \item $\eta:I\to A$ is a morphism in $\CC$, called \emph{unit}, such that the following diagram commutes:
    
    \[\xymatrix{
    A\otimes I \ar[r]^{\id_A\otimes\eta} \ar[dr]_{\id_A} & A\otimes A \ar[d]_\mu & I\otimes A \ar[l]_{\eta\otimes\id_A} \ar[dl]^{\id_A}\\
     & A. &
    }\]
This property is called \emph{unitality} of the product.
\end{itemize}
In several texts in the literature, for instance in \cite{mpsw}, monoids are called \emph{algebras}.
\end{defi}

\begin{defi}
A \emph{comonoid} in $\CC$ is a monoid in $\CC^\op$.
\end{defi}

\begin{rema}
A monoid in a monoidal category is a monoid in its transpose, and viceversa. The same happens for comonoids.
\end{rema}

\begin{defi}
Let $(\CC,\otimes,I)$ be a monoidal category, and let $(A,\mu,\eta)$ and $(A',\mu',\eta')$ be monoids in $\CC$. A \emph{morphism of monoids from $A$ to $A'$} is a morphism $f:A\to A'$, such that $f\circ\mu= \mu'\circ(f\otimes f)$ and $f\circ\eta=\eta'$.
\end{defi}

\begin{rema}
It can be proved directly that identities and composition of morphisms of monoids are morphisms of monoids. Thus, they form a category, often named $\text{Mon}(\CC)$.
\end{rema}

Next we will give examples of monoidal categories which give rise to (co)monoids which are usually found in the literature.

\begin{ex}
For the monoidal category of $\Bbbk$-vector spaces $(\Vect, \otimes_\Bbbk, \Bbbk)$, (co)monoids are ordinary $\Bbbk$-(co)algebras.
\end{ex}

\begin{ex}
Let $R$ be a $\Bbbk$-algebra. There is a monoidal structure in the category of $R$-bimodules $(_R\Mod_R, \otimes_R, R)$, where $R$ acts on itself by multiplication. Its monoids are called \emph{$R$-algebras}.
\end{ex}

\begin{ex}\label{exVectZ}
The category of $\Z$-graded $\Bbbk$-vector spaces has a monoidal structure $(\Vect^\Z, \otimes, \Bbbk_0)$, where the product of two $\Z$-graded vector spaces $V=\bigoplus_{n\in\Z}V_n$ and $W=\bigoplus_{n\in\Z}W_n$ is given by
\[
(V\otimes W)_n = \bigoplus_{p+q=n}V_p\otimes_\Bbbk W_q,
\]
and the unit $\Bbbk_0$ is the graded vector space $\Bbbk$ concentrated in degree $0$. (Co)monoids in this monoidal category are known as \emph{graded (co)algebras}.
\end{ex}

\begin{ex}\label{defitensorCoch}
 Let $(C^\bullet, d_C^\bullet)$ and $(D^\bullet, d_D^\bullet)$ be two
 nonnegatively graded cochain complexes of vector spaces. The product $C\otimes D$ of their underlying $\N_0$-graded spaces can be defined as in the previous example, and the map $\delta^\bullet$ defined by
\[
\delta^{p+q}|_{C^p\otimes D^q} = d_C^p\otimes\id_{D^q} + (-1)^p\ \id_{C^p}\otimes d_D^q
\]
is a differential.
This happens to be the total complex of the product of both complexes. A unit element for this product is again the graded vector space $\Bbbk_0$ with zero differential, giving rise to the monoidal category $(\Coch(\Vect),\otimes,\Bbbk_0)$. (Co)monoids in this category are called \emph{differential graded (co)algebras}.

One can define a monoidal structure in the category of chain complexes in an analogous way.
\end{ex}

\begin{rema}\label{remaChCochCC}
The previous constructions can be made for $\CC^\Z, \Ch(\CC), \Coch(\CC)$ whenever $\CC$ is a monoidal category with a $\Bbbk$-linear structure and arbitrary direct sums.
\end{rema}

\begin{ex}
Let $(H,\mu,\eta,\Delta,\varepsilon)$ be a $\Bbbk$-bialgebra. It is well-known that the category $_H\Mod$ of left $H$-modules is a monoidal category with the usual tensor product of vector spaces. 

(Co)-monoids in this monoidal category are called \emph{$H$-module (co)algebras}.
\end{ex}

\begin{ex}\label{exHcoMod}
Let $(H,\mu,\eta,\Delta,\varepsilon)$ be a $\Bbbk$-bialgebra. The category $^H\Mod$ of left $H$-comodules is a monoidal category with the usual tensor product of vector spaces.

(Co)-monoids in this monoidal category are called \emph{$H$-comodule (co)algebras}.

If $H = \Bbbk G$ is a group algebra, the equality
\[
V_g=\{v\in V:\ v_{[-1]}\otimes v_{[0]}=g\otimes v \}
\]
gives a $1-1$ correspondence between left $H$-comodules and \emph{$G$-graded vector spaces}, that is, vector spaces with a decomposition $V=\bigoplus_{g\in G}V_g$. The formula above for the tensor product of two $H$-comodules $V$ and $W$ yields a grading in the tensor product $V\otimes W=\bigoplus_{g\in G}(V\otimes W)_g$, where $(V\otimes W)_g = \bigoplus_{st=g}V_s\otimes W_t$. The grading on the unit $\Bbbk$ induced by its $\Bbbk G$-coaction is concentrated in degree $1_G$.
In this case, an $H$-comodule algebra is known as a \emph{$G$-graded algebra}, that is, a $\Bbbk$-algebra $A=\bigoplus_{g\in G} A_g$ such that $A_s A_t\subseteq A_{st},\ \forall s,t\in G$, and $1_A\in A_{1_G}$.
If $G=\Z$, the monoidal category obtained this way is exactly the category $\Vect^\Z$ of Example \ref{exVectZ}. Since this case includes the examples we are focused on, in order to avoid confusion with homological degree, we will say that an element $v$ of a $G$-graded vector space $V$ has \emph{internal degree} $g$ if $v\in V_g$.

\end{ex}

\begin{ex}\label{exHYD}
Let $(H,\mu,\eta,\Delta,\varepsilon, \ant)$ be a $\Bbbk$-Hopf algebra. A \emph{left-left Yetter-Drinfeld module over H} is a vector space $V$ endowed with a left $H$-action and a left $H$-coaction such that
\[
(h\cdot v)_{[-1]} \otimes (h\cdot v)_{[0]} = h_{(1)}v_{[-1]}\ant (h_{(3)})\otimes h_{(2)}\cdot v_{[0]},\ \forall h\in H, v\in V.
\]
A \emph{morphism of left-left Yetter-Drinfeld modules over $H$} is a linear map that is simultaneously a morphism of left $H$-modules and a morphism of left $H$-comodules.
The category consisting of the class of left-left Yetter Drinfeld modules over $H$ and their morphisms is denoted by $\HYD$. It can be given a monoidal structure with the usual tensor of vector spaces, defining an action and a coaction on the tensor product of two Yetter-Drinfeld modules in the same way as in the previous examples. It is straightforward to check that the new action and coaction are compatible.

If $H= \Bbbk G$ is a group algebra, the compatibiliy condition translates as $g\cdot V_h \subseteq V_{ghg^{-1}}$. If, in addition, $G$ is abelian, this condition results in $g\cdot V_h\subseteq V_h$, so a Yetter-Drinfeld module over $\Bbbk G$ is a $G$-module $V=\bigoplus_{h\in G}V_h$, where $V_h$ is a $G$-submodule of $V$ for every $h\in G$.
\end{ex}

\begin{ex}\label{exJP}
The \emph{Jordan plane} is the algebra $A=\frac{\Bbbk\langle x,y\rangle}{\left<yx-xy+\frac{1}{2}\ x^2\right>}$. 

Let us define an action of the Hopf algebra $\Bbbk\Z=\Bbbk[t,t^{-1}]$ on $A$ by $t\cdot x= x,\ t\cdot y=x+y$, and a $\Z$-grading on $A$ setting $x,y\in A_1$. That can be extended multiplicatively to $\Bbbk\langle x,y\rangle$, and are well defined after taking the quotient because the relation is homogeneous and preserved by the action. Since the action preserves the grading, $A$ is a Yetter-Drinfeld module over $\Bbbk\Z$. By definition, the algebra structure of $A$ is compatible with the action and the coaction of $\Bbbk\Z$, so $A$ is a monoid in the monoidal category $\kZYD$.
\end{ex}

\begin{ex}\label{exSJP}
This is an example of the same type, with a different action of $\Z$. The \emph{super Jordan plane} is the algebra $A=\frac{\Bbbk\langle x,y\rangle}{\left<x^2, y^2 x - xy^2 - xyx\right>}$ with $\Bbbk\Z$-action on $A$ given by $t\cdot x= -x,\ t\cdot y=x-y$, and $\Z$-grading on $A$ setting $x,y\in A_1$. As in the previous example, they induce on $A$ the structure of Yetter-Drinfeld module over $\Bbbk\Z$ which is compatible with its product and unit, so $A$ is a monoid in $\kZYD$.
\end{ex}

\begin{ex}\label{exomega}
Let $\rho:A\to R$ be a $\Bbbk$-algebra homomorphism, and let $P_\bullet\xrightarrow{f} A$ be a projective resolution of $A$ as $A$-bimodule. The triple $(P_\bullet,\omega,f)$ is a comonoid in the monoidal category $\left(\overline{\Ch}\left(\AModA\right),\otimes_A, A_0\right)$, where the bar indicates that morphisms in this category are those of $\Ch\left(\AModA\right)$ modulo homotopy (this is necessary for the associativity and unitality conditions).

The convolution product makes $\bigl(\Hom_{AA}(P_\bullet,R),\cup,(1\mapsto\rho\circ f)\bigr)$ a monoid in the monoidal category $(\overline{\Coch}(\Vect),\otimes,\Bbbk_0)$, and taking cohomology makes $\H^\bullet(A,R)$ a monoid in the monoidal category $(\Vect^\Z,\otimes, \Bbbk_0)$.
\end{ex}

\begin{defi}
Let $(\CC,\otimes,I)$ and $(\DD,\times,J)$ be monoidal categories. A \emph{lax monoidal functor from $\CC$ to $\DD$} is a triple $(F,\varphi,\varphi_0)$, where

\begin{itemize}
    \item $F:\CC\to\DD$ is a functor,
    \item $\varphi$ is a natural transformation given by a family of morphisms $\{\varphi_{X,Y}\}_{X,Y\in\CC}$, where $\varphi_{X,Y}:F(X)\times F(Y)\to F(X\otimes Y)$, such that for every $X,Y,Z\in\CC$ the following diagram commutes:
\[
\xymatrix{
F(X)\times F(Y)\times F(Z)\ar[rr]^-{\varphi_{X,Y}\times \id_Z} \ar[d]_{\id_X\times\varphi_{Y,Z}} & & F(X\otimes Y)\times Z \ar[d]^{\varphi_{X\otimes Y, Z}} \\
F(X)\times F(Y\otimes Z) \ar[rr]_{\varphi_{X,Y\otimes Z}} & & F(X\otimes Y\otimes Z),
}
\]
    \item $\varphi_0:J\to F(I)$ is a morphism in $\DD$ such that the following diagrams commute:
    
\[
\xymatrix{
F(X) \ar[r]^-= \ar[d]_= & J\times F(X) \ar[d]^{\varphi_0\times \id_{F(X)}}\\
 F(I\otimes X)  & F(I)\times F(X), \ar[l]^{\varphi_{I,X}}
} \quad\quad
\xymatrix{
F(X) \ar[r]^-= \ar[d]_= &  F(X)\times J \ar[d]^{\id_{F(X)}\times\varphi_0}\\
 F(X\otimes I)  & F(X)\times F(I). \ar[l]^{\varphi_{X,I}}
 }
\]
We will denote a lax monoidal functor by $(F,\varphi, \varphi_0):(\CC, \otimes, I)\to (\DD, \times, J)$, and eventually call it a \emph{lax monoidal structure on F}.
\end{itemize}
\end{defi}

\begin{defi}
A \emph{strong monoidal functor} is a lax monoidal functor $(F,\varphi, \varphi_0)$ such that $\varphi$ is a natural isomorphism and $\varphi_0$ is an isomorphism.
\end{defi}

\begin{ex}
We may consider the functor $V\mapsto V^*=\Homk(V,\Bbbk)$ as a monoidal functor $((-)^*, \varphi, \varphi_0): (\Vect^\op, \otimes, \Bbbk) \to (\Vect, \otimes, \Bbbk)$. The natural transformation $\varphi$  given by $\varphi_{V,W}:V^*\otimes W^* \to (V\otimes W)^*$ is defined by $\varphi_{V,W}(f\otimes g) = f\otimes g$, considering the first $f\otimes g$  as the tensor product of two elements and the last one as the tensor product of two morphisms. The morphism $\varphi_0=:\Bbbk\to\Bbbk^*$ is defined by $\varphi_0(1)=\id_\Bbbk$. This monoidal functor is not strong, but its restriction to the opposite category of finite dimensional vector spaces is so.
\end{ex}

\begin{rema}\label{remaforgetstrong}
Defining monoidal structures like in Examples \ref{exVectZ} to \ref{exHYD}, adding extra structure to the usual tensor product of $\Bbbk$-vector spaces, automatically gives their respective forgetful functors strong monoidal structures via identities.
\end{rema}

\begin{ex}
In the case considered in Example \ref{exHYD} where $H=\Bbbk\Z$, the correspondence between left-left Yetter Drinfeld modules over $\Bbbk\Z$ and $\Z$-graded $\Z$-modules can be stated by saying that the identity induces a strong monoidal functor from $\kZYD$ to $\cc_{\Bbbk\Z}\Mod\dd^\Z$, where the monoidal structure of the latter is the one given in Remark \ref{remaChCochCC}.
\end{ex}

\begin{ex}
Let $R$ be a $\Bbbk$-algebra. The forgetful functor $U:\,_R\Mod_R\to\Vect$ admits a lax monoidal structure $(U,\pi,\eta):\left(_R\Mod_R,\otimes_R, R\right)\to (\Vect,\otimes,\Bbbk)$, where $\pi_{M,N}:M\otimes N\to M\otimes_R N$ is the canonical projection and $\eta:\Bbbk\to R$ is the unit morphism of $R$. This lax monoidal functor is not strong in general.
\end{ex}

\begin{ex}
The K\"unneth Theorem for chain complexes of vector spaces states that the homology functor $H^\bullet:\Ch(\Vect)\to\Vect^\Z$ is strong monoidal.
\end{ex}

\subsection{Braided monoidal categories}

In this section we will pay attention to an additional structure in some monoidal categories: braidings.

\begin{defi}
A \emph{braided monoidal category} is a 7-tuple $(\CC,\otimes,I, a,\ell,r,c)$, where

\begin{itemize}
    \item $(\CC, \otimes, I, a, \ell, r)$ is a monoidal category,
    \item $c$ is a natural isomorphism, called \emph{braiding}, given by a family of morphisms $c_{X,Y}: X\otimes Y \to Y\otimes X$ for each $X,Y\in\CC$, such that the following diagrams commute for every $X,Y,Z\in\CC$:
\[
\xymatrix{X\otimes (Y\otimes Z) \ar[rr]^{c_{X, Y\otimes Z}} \ar[d]_{a_{X,Y,Z}^{-1}}  && (Y\otimes Z)\otimes X \\
(X\otimes Y)\otimes Z \ar[d]_{c_{X,Y}\otimes\id_Z} & & Y\otimes (Z\otimes X) \ar[u]_{a_{Y,Z,X}^{-1}}\\
(Y\otimes X)\otimes Z \ar[rr]_{a_{Y,X,Z}}  && Y\otimes (X\otimes Z), \ar[u]_{\id_Y\otimes c_{X,Z}}}
\]

\[
\xymatrix{(X\otimes Y) \otimes Z \ar[rr]^{c_{X\otimes Y, Z}} \ar[d]_{a_{X,Y,Z}}  && Z\otimes (X\otimes Y) \\
X\otimes(Y\otimes Z) \ar[d]_{\id_X\otimes c_{Y,Z}} && (Z\otimes X) \otimes Y \ar[u]_{a_{Z,X,Y}}\\
X\otimes (Z\otimes Y) \ar[rr]_{a_{X,Z,Y}^{-1}}  && (X\otimes Z) \otimes Y. \ar[u]_{c_ {X,Z}\otimes\id_Y}}
\]
\end{itemize}
\end{defi}

\begin{rema}
It is proven in Proposition 1.1 of \cite{joyal1993braided} that commutativity of the previous diagrams implies commutativity of the following diagram for every $X\in\CC$:
\[
\xymatrix{
X\otimes I \ar[rr]^{c_{X,I}} \ar[dr]_{r_X} & & I\otimes X \ar[dl]^{\ell_X} \\ &X&}.
\]
\end{rema}

\begin{rema}
Let $(\CC,\otimes,I,a,\ell,r,c)$ be a braided monoidal category, and consider the monoidal structures on its opposite and transpose categories as in Remark \ref{remaflop} and Definition \ref{defitranspose}, respectively. They both admit the braiding $\Tilde{c}$, where $\Tilde{c}_{X,Y}=c_{Y,X}$ for all $X,Y\in\CC$.
\end{rema}

\begin{rema}
For a braided monoidal category $(\CC, \otimes, I, a, \ell, r, c)$, inverting the braiding gives another braided monoidal category $(\CC, \otimes, I, a, \ell, r, c^{inv})$, where $c^{inv}_{X,Y}=\cc c_{Y,X}\dd\inv$ for all $X,Y\in\CC$.
\end{rema}

\begin{rema}\label{remabraideqn}
If a braided monoidal category is strict as a monoidal category, the hexagonal diagrams above turn into triangles, yielding the following equalities for every $X,Y,Z\in\CC$:
\[
c_{X,Y\otimes Z} = \cc \id_Y\otimes c_{X,Z}\dd\circ\cc c_{X,Y}\otimes\id_Z\dd,
\] \[
c_{X\otimes Y,Z} = \cc c_{X,Z}\otimes\id_Y\dd \circ\cc\id_X\otimes c_{Y,Z}\dd.
\]
From these equalities one can deduce the \emph{braid equation}:
\[
\cc c_{Y,Z}\otimes\id_X\dd \circ \cc\id_Y\otimes c_{X,Z}\dd \circ \cc c_{X,Y}\otimes\id_Z\dd =
\cc\id_Z\otimes c_{X,Y}\dd \circ \cc c_{X,Z}\otimes\id_Y\dd \circ \cc\id_X\otimes c_{Y,Z}\dd.
\]
All braided monoidal categories considered in this work will be strict as monoidal categories.
\end{rema}

\begin{defi}
Let $(\CC,\otimes,I,c)$ be a braided monoidal category.
\begin{itemize}
    \item A monoid $(A,\mu,\eta)$ in $\CC$ is \emph{commutative} if $\mu\circ c_{A,A}=\mu$.
    \item A comonoid $(A,\Delta,\varepsilon)$ in $\CC$ is \emph{cocommutative} if $c_{A,A}\circ\Delta=\Delta$.
\end{itemize}
\end{defi}

Next we will discuss some relevant examples.

\begin{ex}\label{extau}
The formula $\tau_{V,W}(v\otimes w)=w\otimes v, \ \forall v\in V, w\in W$, defines a braiding $\tau$ in the categories $\Vect$, $^H\Mod$ for a commutative $\Bbbk$-bialgebra $H$, and $_H\Mod$ for a cocommutative $\Bbbk$-bialgebra $H$. In these cases, (co)commutative (co)monoids are monoids in their respective monoidal categories whose underlying $\Bbbk$-(co)algebra is (co)commutative in the usual sense.
\end{ex} 

\begin{rema}\label{remabraidedcomm}
When working in braided monoidal categories whose objects are $\Bbbk$-vector spaces with additional structure, commutativity may be called \emph{braided commutativity}, and analogously with cocommutativity. This is useful to emphasize that the braiding involved is not $\tau$, and so the notion of commutativity involved differs from the usual one.
\end{rema}

\begin{ex}\label{extaugr}
In the categories $\Ch(\Vect), \Coch(\Vect)$ and $\Vect^\Z$, one can define a braiding $\tau^{gr}$ by the formula $\tau^{gr}(v\otimes w)=(-1)^{pq}\ w\otimes v, \ \forall v\in V_p,w\in W_q$. (Co)commutativity for this braiding is called \emph{graded (co)commutativity}. In the categories $\hCh(\Vect)$ and $\hCoch(\Vect)$, it is called \emph{graded (co)commutativity up to homotopy}.
\end{ex}

\begin{ex}
The previous construction can be generalized to an arbitrary braided monoidal category $(\CC,\otimes,I,c)$ with a $\Bbbk$-linear structure: in this case, one can define a braiding $c^{gr}$ in the categories $\Ch(\CC),\Coch(\CC)$ and $\CC^\Z$ by the formula $c^{gr}_{V,W}|_{V_p\otimes W_q} = (-1)^{pq}c_{V_p,W,q}$. As in Remark \ref{remabraidedcomm}, (co)commutativity in these categories may be called \emph{graded braided (co)commutativity}, or \emph{braided graded (co)commutativity}, and also its versions up to homotopy in $\hCh(\CC)$ and $\hCoch(\CC)$.
\end{ex}

\begin{ex}\label{exbraidHYD}
Let $H$ be a $\Bbbk$-Hopf algebra. The category $(\HYD,\otimes,\Bbbk)$ admits a braiding $c$ such that
\[
c_{V,W}(v\otimes w) = v_{[0]}\cdot w\otimes v_{[1]},\ \forall v\in V,w\in W.
\]
In the case $H=\Bbbk G$ for a group $G$, the braiding has the following expression:
\[
c_{V,W}(v\otimes w) = g\cdot w\otimes v,\ \forall g\in G, v\in V_g,w\in W.
\]
\end{ex}

In the following sections we will define certain (co)monoids in the monoidal categories $\hCh\cc\HYD\dd, \hCoch\cc\HYD\dd$ and $\cc\HYD\dd^\Z$, and prove that they are (co)commutative. To introduce them, we will need the following concept, which generalizes the definition of bialgebra to arbitrary braided monoidal categories.

\begin{defi}
Let $(\CC,\otimes,I,c)$ be a braided monoidal category. A \emph{bimonoid} in $\CC$ is a 5-tuple $(A,\mu,\eta,\Delta,\varepsilon)$, where
\begin{itemize}
    \item $(A,\mu,\eta)$ is a monoid,
    \item $(A,\Delta, \varepsilon)$ is a comonoid,
    \item the following diagrams commute:
\[
\xymatrix{
A\otimes A \ar[r]^-\mu \ar[d]_{\Delta\otimes\Delta} & A \ar[r]^-\Delta & A\otimes A \\ A\otimes A\otimes A\otimes A \ar[rr]_{\id_A\otimes c_{A,A}\otimes\id_A}  & & A\otimes A\otimes A\otimes A,\ar[u]_{\mu\otimes\mu}}
\]
\[
\xymatrix{
I \ar[r]^-{\simeq} \ar[d]_{\eta} & I\otimes I \ar[d]^{\eta\otimes\eta} \\
A \ar[r]_-{\Delta} & A\otimes A,
} \quad 
\xymatrix{
A\otimes A \ar[r]^-{\mu} \ar[d]_\varepsilon & A \ar[d]^{\varepsilon\otimes\varepsilon} \\
I\otimes I \ar[r]_{\simeq}  & I,
} \quad
\xymatrix{
I \ar[rr]^{\id_I} \ar[dr]_\eta & & I \\ & A.\ar[ur]_\varepsilon &
}
\]
\end{itemize}
\end{defi}

\begin{rema}\label{remabimonoideqv}
In a braided monoidal category $(\CC,\otimes,I,c)$ for any two monoids $(A,\mu,\eta)$ and $(A',\mu',\eta')$, the triple $(A\otimes A', (\mu\otimes\mu')\circ(\id_A\otimes c_{A',A}\otimes\id_{A'}), \eta\otimes\eta')$ is also a monoid.

Analogously, if $(A,\Delta, \varepsilon)$ and $(A',\Delta', \varepsilon')$ are two comonoids, their product admits the comonoid structure $(A\otimes A', (\id\otimes c_{A,A'}\otimes\id_{A'})\circ (\Delta\otimes\Delta'), \varepsilon\otimes\varepsilon')$.

As for ordinary bialgebras, the axioms of a bimonoid $(A,\mu,\eta,\Delta,\varepsilon)$ are equivalent to $\Delta$ and $\varepsilon$ being morphisms of monoids, and they are also equivalent to $\mu$ and $\eta$ being morphisms of comonoids.
\end{rema}

\begin{defi}
A \emph{Hopf monoid} is a 6-tuple $(A,\mu,\eta,\Delta,\varepsilon,s)$, where $(A,\mu,\eta,\Delta,\varepsilon)$ is a bimonoid, and $s:A\to A$ is a morphism in $\CC$, called \emph{antipode}, such that
\[
\mu\circ\cc\id_A\otimes s\dd\circ\Delta = \eta\circ\varepsilon = \mu\circ\cc s\otimes\id_A\dd\circ\Delta.
\]
\end{defi}

\begin{rema}
Bimonoids and Hopf monoids in a braided monoidal category $\CC$ are often named \emph{braided bialgebras} and \emph{braided Hopf algebras}, respectively.
\end{rema}

\begin{rema}
Unlike previous examples, the so called graded Hopf algebras in the literature do not correspond to Hopf monoids in the category $\Vect^\Z$, because the braiding used in their definition is just the braiding of vector spaces.
\end{rema}
\bigskip

\subsection{Nichols Algebras}
In this section we will recall the construction of the Nichols algebra of a Yetter-Drinfeld module $V$ over a Hopf algebra $H$. Nichols algebras are examples of braided Hopf algebras in the category $\HYD$ for some Hopf algebra $H$.

Let $T(V)=\bigoplus_{n\in\N_0}V^{\otimes n}$ denote the tensor algebra of $V$ with the concatenation product and the inclusion of $\Bbbk=V^{\otimes 0}$ as unit. As seen in Example \ref{exHYD}, $T(V)$ can be regarded as a Yetter-Drinfeld module over $H$, and it is immediate to check that the aforementioned product and unit preserve the Yetter-Drinfeld module structure. Therefore, one can consider $T(V)$ as a monoid in $\HYD$.

Now, the braiding of this category induces a monoid structure on $T(V)\otimes T(V)$, which by Remark \ref{remaforgetstrong} can be seen as an ordinary algebra structure.

Since $T(V)$ is the free algebra generated by $V$, defining $\Delta(v)=1\otimes v + v\otimes 1$ extends to a unique algebra morphism $\Delta:T(V)\to T(V)\otimes T(V)$. As a consequence of the definition of the $H$-action and $H$-coaction on tensor products, the map $\Delta$ turns out to be a morphism of Yetter-Drinfeld modules. The map $\varepsilon:T(V)\to\Bbbk$ defined as the projection onto $V^{\otimes 0}$ is a counit for $\Delta$, and trivially preserves the algebra structure and the Yetter-Drinfeld module structure. Therefore, $T(V)$ with this structure is a bimonoid in $\HYD$, which happens to be a Hopf monoid.

There exists a maximal Hopf ideal and Yetter-Drinfeld submodule $J$ of $T(V)$ among those contained in $\bigoplus_{n\geq 2}V^{\otimes n}$, see \cite[Proposition 2.2]{andruskiewitsch2002pointed}. The \emph{Nichols algebra of $V$} is the quotient $\BB(V)=\frac{T(V)}{J}$.

The Jordan plane of Example \ref{exJP} and the super Jordan plane of Example \ref{exSJP} are the Nichols algebras corresponding to the $\Bbbk\Z$-Yetter Drinfeld modules, called $\mathcal{V}(1,2)$ and $\mathcal{V}(-1,2)$ respectively.

The Nichols algebra $\BB(V)$ also inherits a natural $\Z$-grading from $T(V)$, which is preserved by its product and coproduct. In the Jordan and super Jordan plane it coincides with the $\Z$-grading given by their $\Bbbk\Z$-comodule structure.

\section{Graded braided (co)commutativity}

The aim of this section is to provide a setting where the Hochschild cohomology of $A$ with coefficients in $\Bbbk$
could be given a graded Yetter-Drinfeld module structure over $H$, where $A$ is a Nichols algebra in the category $\HYD$. In this context, the $\Z$-graded vector space $\H^\bullet(A,\Bbbk)$ can be seen as a monoid in the category $\cc\HYD\dd^\Z$, which happens to be graded braided commutative. We will follow the notation in Section 3 of \cite{mpsw}, whose main result is recovered in Corollary \ref{corMPSW}.

\begin{defi}\label{defihomint}
Let $(\CC, \otimes, I)$ be a monoidal category. Let $V,Y$ be objects in $\CC$. An \emph{internal hom} is an object in  $\CC$, denoted $\hom(V,Y)$, equipped with a natural isomorphism
\[
\Theta_{V,Y}:\Hom_\CC(-\otimes V, Y) \longrightarrow \Hom_\CC(-, \hom(V,Y)).
\]
For $X\in\CC$, we will write $\Theta_{X,V,Y} := \cc\Theta_{V,Y}\dd_X$.
\end{defi}

\begin{rema}\label{remaYoneda}
Wherever $\hom(V,Y)$ exists, it is unique up to isomorphism and it is functorial in both variables, due to Yoneda's Lemma. See, for example, \cite{maccategories}, Chapter III, Section 2. As a consequence, if $V\in\CC$ is such that $\hom(V,Y)$ exists for all $Y\in\CC$, then $\hom(V,-)$ is a functor which is right adjoint to $-\otimes V$. 
\end{rema}
\begin{defi}
If the adjunction $\Theta_{V,-}$ of the previous definition exists for all $V\in\CC$, the monoidal category $(\CC,\otimes,I)$ is called a \emph{right closed monoidal category}.
\end{defi}  

\begin{ex}\label{exhomintVectZ}
The monoidal category $(\Vect^\Z,\otimes,\Bbbk_0)$ is right closed monoidal, its internal hom functor is called \emph{graded hom}. For a graded $\Bbbk$-vector space $Y$ and $n\in\Z$, let $Y[n]$ denote the graded vector space such that $Y[n]_i=Y_{i+n}$. The graded hom object $\hom(V,Y) = \bigoplus_{n\in\Z}\hom(V,Y)_n$, is defined as follows:
\[
\hom(V,Y)_n = \Hom_{\Vect^\Z}(V,Y[n]) = \prod_{i\in\Z}\Homk(V_i,Y_{i+n}),\quad\forall\ V,Y\in\Vect^\Z.
\]
\end{ex}

\begin{rema}\label{remahomgrincHom}
The canonical inclusion $\hom(V,Y)\subseteq\Homk(V,Y)$ is an equality if and only if $V$ is finite dimensional. This is equivalent to the condition that $V_n$ is finite dimensional for all $n$, and nonzero only for finitely many values of $n$.
\end{rema}

\bigskip

The internal hom objects in the following examples can be defined more generally for a Hopf algebra $H$ with bijective antipode. We will focus in the case $H=\Bbbk G$.

\begin{ex}
The category $(_{\Bbbk G}\Mod, \otimes, \Bbbk)$ is right closed monoidal, with its internal hom defined by $\hom(V,Y)=\Homk(V,Y)$ and $G$-action given by
\[
(g\cdot f)(x) = g\cdot f(g^{-1}\cdot x),\ \forall g\in G,\ \forall x\in V,\ \forall f\in \Homk(V,Y)
\]
\end{ex}

\begin{ex}\label{exhomintVectG}
Let $G$ be a group. The construction for $\Z$-graded vector spaces in Example \ref{exhomintVectZ} can be generalized to the monoidal category $(^{\Bbbk G}\Mod, \otimes, \Bbbk)$, which is isomorphic to the category $\Vect^G$ of $G$-graded vector spaces -- see Example \ref{exHcoMod}. It has an internal hom functor defined by $\hom(V,Y) = \bigoplus_{h\in G}\hom(V,Y)_h$, where
\[
\hom(V,Y)_h = \prod_{s\in G} \Homk(V_s, Y_{hs}).
\]
As in Remark \ref{remahomgrincHom}, there is a canonical inclusion $\hom(V,Y)\subseteq\Homk(V,Y)$, which is an equality either when $V$ is finite dimensional or when $G$ is finite.
\end{ex}

\begin{ex} \label{exhomintkGYD}
The category $(\kGYD, \otimes, \Bbbk)$ of Yetter-Drinfeld modules over a group algebra has an internal hom which combines the $\Bbbk G$-module and $\Bbbk G$-comodule structures of its objects. For $V,Y\in\kGYD$, the $G$-graded space $\hom(V,Y)$ is defined as in the previous example, and the action of $g\in G$ on $f=(f_s)_{s\in G}\in\hom(V,Y)_h$ is given by
\[
(g\cdot f)_s(v) = g\cdot (f_{g^{-1}sg}(g^{-1}\cdot v)),\ \forall s\in G, \ \forall v\in V_s.
\]
It is straightforward to check that $g\cdot f\in\hom(V,Y)_{ghg\inv} \ \forall g,h\in G, f\in\hom(V,Y)_h$, so $\hom(V,Y)$ is indeed an object in $\kGYD$.
It also holds that $\hom(V,Y)\simeq\Homk(V,Y)$ when $V$ is finite dimensional or $G$ is finite.
\end{ex}

 The following remark is generalized in \cite{mpsw} to arbitrary monoidal categories. Their aim is, when possible, to replace the functor $\Homk(-,\Bbbk)$ by an internal hom functor, obtaining the desired construction in the cases in which they coincide.

\begin{rema}
Let $(A,\varepsilon)$ be an augmented $\Bbbk$-algebra, and consider the cochain complex $(C^\bullet(A), d^\bullet)$. Taking into account the action of $A$ on $\Bbbk$, the coboundary $d^n:C^n(A)\to C^{n+1}(A)$ is as follows:
\[
\begin{array}{l}

d^n f(a_1\otimes\cdots\otimes a_{n+1})\ =\ \varepsilon(a_1) f(a_2\otimes\cdots\otimes a_{n+1})\\
\quad\quad +\ \sum_{i=1}^n (-1)^i\ f(a_1\otimes\cdots\otimes a_i a_{i+1}\otimes\cdots\otimes a_{n+1}) + (-1)^{n+1}\ f(a_1 \otimes\cdots\otimes a_n)\varepsilon(a_{n+1})\\ 
\quad =\ f \big{(}\ \varepsilon(a_1) a_2\otimes\cdots\otimes a_{n+1} + \sum_{i=1}^n (-1)^i\ a_1\otimes\cdots\otimes a_i a_{i+1}\otimes\cdots\otimes a_{n+1}\\
\quad\quad +\ (-1)^{n+1}\ a_1 \otimes\cdots\otimes a_n\varepsilon(a_{n+1})\ \big{)}.
\end{array}
\]
This can also be read as $d^n f=f\circ\delta_{n}$, where $\delta_n: A^{\otimes n+1}\to A^{\otimes n}$ is 
\[
\begin{array}{l}
\delta_n(a_1\opo a_{n+1})\ =\  \varepsilon(a_1) a_2\otimes\cdots\otimes a_{n+1}\\
\qquad +\ \sum_{i=1}^n (-1)^i a_1\otimes\cdots\otimes a_i a_{i+1}\otimes\cdots\otimes a_{n+1}\ +\ (-1)^{n+1} a_1 \otimes\cdots\otimes a_n\varepsilon(a_{n+1}).
\end{array}
\]
Setting $S_n(A)=A^{\otimes n}$ makes $(S_\bullet(A),\delta_\bullet)$ a chain complex of $\Bbbk$-vector spaces satisfying the equality $(C^\bullet(A),d^\bullet)=(\Homk(S_\bullet(A),\Bbbk), (\delta_\bullet)^*)$.
\end{rema}

Next we follow the general construction given in \cite[Section 3]{mpsw}.

\begin{defi}
Let $(\CC,\otimes,I)$ be a monoidal category, where $I$ is a monoid with the canonical isomorphism $I\otimes I\simeq I$. An \emph{augmented monoid in $\CC$} is a 4-tuple $(A,\mu,\eta,\varepsilon)$, where $(A,\mu,\eta)$ is a monoid in $\CC$ and $\varepsilon:A\to I$ is a morphism of monoids. As usual, we omit part of the structure when it is clear from the context.
\end{defi}

\begin{rema}
If $(\CC,\otimes,I,c)$ is a braided monoidal category and $(A,\mu,\eta,\Delta,\varepsilon)$ is a bimonoid in $\CC$, then $(A,\mu,\eta,\varepsilon)$ is an augmented monoid in $\CC$.
\end{rema}

\begin{defi}\label{defiS}
Let $(A,\mu,\eta,\varepsilon)$ be an augmented monoid in a monoidal category $(\CC,\otimes,I)$. The triple $(S(A),\{\partial_i\}_i,\{\sigma_j\}_j)$ defined as follows is a simplicial object in $\CC$:

\begin{itemize}
    \item Objects: $S_n(A) = A^{\otimes n}$.
    \item Faces: $\partial_i\colon S_{n+1}(A)\to S_n(A)$ 
    
    \[\partial_i = \left\{ \begin{array}{lr}
       \varepsilon\otimes\id_A^{\otimes n}  & \text{if } i=0 \\
        \id_A^{\otimes i-1}\otimes\mu\otimes\id_A^{n-i}\quad & \text{if } i=1,\ldots,n \\
        \id_A^{\otimes n}\otimes\varepsilon & \text{if } i=n+1.
    \end{array}
    \right.\]
    
     \item Degeneracies: $\sigma_j\colon S_n(A)\to S_{n+1}(A)$
     
     \[
     \sigma_j = \id_A^{\otimes j}\otimes\eta\otimes\id_A^{\otimes n-j},\quad \forall j=0,\ldots,n.
     \]
\end{itemize}

We will denote its associated chain complex in $\CC$ by $(S_\bullet(A),\delta_\bullet)$, where the differential $\delta_{n}:S_{n+1}(A)\to S_n(A)$ is defined as $\delta_n=\sum_{i=o}^{n+1}(-1)^i\partial_i$.
\end{defi}

\begin{rema}
This construction is functorial, therefore it gives functors $S$ and $S_\bullet$ from the category of augmented monoids in $\CC$ to the categories of simplicial objects in $\CC$ and of chain complexes in $\CC$, respectively.
\end{rema}

\begin{defi} \cite[Definition 3.4]{mpsw}
Let $A$ be an augmented monoid in a monoidal category $(\CC,\otimes,I)$. The \emph{Hochschild cohomology of A with trivial coefficients} is the cohomology of the complex $(\hom(S_\bullet(A),I),\hom(\delta_\bullet,I))$. It is denoted by $\H^\bullet(A)$.
\end{defi}

The following are the final steps to obtain a monoid structure in $\CC^\Z$ analogous to the Hochschild cohomology with trivial coefficients on a bialgebra, with the cup product.

\begin{defi}
Let $(\CC,\otimes,I)$ be a monoidal category, and let $V,Y\in\CC$ be such that $\hom(V,Y)$ exists. Let $\Theta$ be the natural isomorphism of Definition \ref{defihomint}. The inverse of the isomorphism
\[
\Theta_{\hom(V,Y),V,Y}:\Hom_\CC(\hom(V,Y)\otimes V,Y)\to\Hom_\CC(\hom(V,Y),\hom(V,Y)).
\]
gives a morphism
\[
ev_{_{V,Y}}=\Theta_{\hom(V,Y),V,Y}^{-1}\cc\id_{\hom(V,Y)}\dd:\hom(V,Y)\otimes V\to Y,
\]
called \emph{evaluation morphism}.
\end{defi}

\begin{defi}\label{defixi} \cite[Remark 3.6]{mpsw}
Let $(\CC,\otimes,I)$ be a monoidal category, and let $V,W\in\CC$ be such that $\hom(V,I), \hom(W,I)$ and $\hom(W\otimes V,I)$ exist. The map
\[
\xi_{V,W}: \hom(V,I)\otimes \hom(W,I) \to \hom(W\otimes V,I)
\]
is defined as follows:
\[
\xi_{V,W} = \Theta_{\hom(V,I)\otimes\hom(W,I),W\otimes V,I} \cc ev_{_{V,I}}\circ(\id_{\hom(V,I)}\otimes ev_{_{W,I}}\otimes\id_{V})\dd.
\]
\end{defi}

\begin{rema}
 If we restrict to those pairs $V,W$ for which $\hom(V,I)\otimes \hom(W,I)$ and $\hom(W\otimes V,I)$ exist, we can say that the collection of maps $\xi_{V,W}$ define a natural transformation $\xi$.
 \end{rema}

\begin{rema}\label{remaxiVect}
As remarked before, in the monoidal categories of Examples \ref{exhomintVectZ} to \ref{exhomintkGYD} the object $\hom(V,Y)$ is naturally included in $\Homk(V,Y)$ for all objects $V$ and $Y$. In these cases, the evaluation morphism turns out to coincide with the usual evaluation of functions. Now, let us consider two elements $f\in\hom(V,\Bbbk)\subseteq\Homk(V,\Bbbk)$ and $g\in\hom(W,\Bbbk)\subseteq\Homk(W,\Bbbk)$. The equality
\[
\xi_{V,W}(f\otimes g) (w\otimes v) = g(w)f(v),\ \forall v\in V, w\in W,
\]
shows that $\xi_{V,W}$ is in fact the canonical inclusion $V^*\otimes W^*\hookrightarrow (W\otimes V)^*$.
\end{rema}

\begin{defi}\label{deficupop}\cite[Definition 3.7]{mpsw}
Let $A$ be an augmented monoid in a monoidal category $(\CC,\otimes,I)$, and let $(S_\bullet(A), \delta_\bullet)$ be the chain complex of Definition \ref{defiS}. The \emph{cup product} $\smallsmile:\H^\bullet(A)\otimes\H^\bullet(A)\to\H^\bullet(A)$ is the collection of morphisms
\[
\smallsmile_{p,q}:\H^p(A)\otimes\H^q(A)\to\H^{q+p}(A)
\]
induced by
\[
\xi_{S_p(A),S_q(A)}:\hom(S_p(A),I)\otimes\hom(S_q(A),I)\to\hom(S_q(A)\otimes S_p(A),I).
\]
\end{defi}

\begin{rema} \label{remacupop}
Let $(\CC,\otimes,\Bbbk)$ be a monoidal category like in Remark \ref{remaxiVect}, and let $A$ be an augmented monoid in $\CC$, which implies that $A$ is also an augmented $\Bbbk$-algebra. If for every $n\in\N_0$ the objects $\hom(A^{\otimes n},\Bbbk)$ exist and their canonical inclusions in $\Homk(A^{\otimes n},\Bbbk)$ are isomorphisms of $\Bbbk$-vector spaces, then the cochain complexes  $(\hom(S_\bullet(A),\Bbbk),\hom(\delta_\bullet,\Bbbk))$ and $(\Homk(S_\bullet(A),\Bbbk),(\delta_\bullet)^*)$ are isomorphic. Taking homology, one gets that $\H^\bullet(A)$ and the usual Hochschild cohomology with trivial coefficients $\H^\bullet(A,\Bbbk)$ are isomorphic as $\Z$-graded $\Bbbk$-vector spaces. In addition, according to the expression of $\xi$, the cup product defined in $\H^{\bullet}(A)$ is the opposite of the usual cup product in $\H^\bullet(A,\Bbbk)$. To avoid confusion, we will use the symbol $\smallsmile$ to refer to the cup product of definition \ref{deficupop}, instead of the symbol $\cup$ used to define the usual cup product in Hochschild cohomology in Section \ref{sectionHochschild}. 
\end{rema}

Next we will recall one of the main results of \cite{mpsw}, together with a corollary. We will prove a generalization of the latter in the end of this section.

\begin{thm} \label{thmMPSW} \cite[Theorem 3.12]{mpsw}
Let $A$ be a bimonoid in a braided monoidal category $(\CC,\otimes,I,c)$, such that $\hom(A^{\otimes n},I)$ exists for all $n\in\N_0$. Then $(\H^\bullet(A),\smallsmile)$ is a graded braided commutative algebra in $\CC^\Z$.
\end{thm}

\begin{cor} \label{corMPSW} \cite[Corollary 3.13]{mpsw}
Let $H$ be a Hopf algebra with bijective antipode, and let $A$ be a bimonoid in $\HYD$. If either $H$ or $A$ is finite dimensional, then $H^\bullet(A,\Bbbk)$ is a $\Z$-graded object in $\HYD$, and the opposite of its cup product is graded braided commutative.
\end{cor}

The authors use this theorem to prove finite generation of the cohomology algebra of some pointed Hopf algebras.

On the other hand, the Jordan and super Jordan plane are infinite dimensional and they are bimonoids in the category of Yetter-Drinfeld modules over $\Bbbk\Z$, which is an infinite dimensional Hopf algebra. The internal hom objects needed in the construction in \cite{mpsw} do not exist, so a priori $\H^\bullet(A,\Bbbk)$ cannot be regarded as a $\Z$-graded Yetter-Drinfeld module over $\Bbbk\Z$, and therefore there is no graded braided commutativity condition to check. However, there is a condition that can be checked, which coincides with graded braided commutativity of the cup product when it can be defined. The next propositions will lead to that condition.

\begin{lemma}
Let $(\CC,\otimes,I,c)$ be a braided monoidal category, and let $V,W\in\CC$. The following diagram commutes:

\[
\scriptsize
\xymatrix{ \hom(V,I)\otimes\hom(W,I)\otimes W\otimes V \ar[rrrr]^{\cc c_{_{\hom(W,I),\hom(V,I)}}\dd\inv\otimes\,c_{_{W,V}}} \ar[dd]_{\id_{_{\hom(V,I)}}\otimes\, ev_{_{W,I}}\otimes\id_{_V}} & & & & \hom(W,I)\otimes\hom(V,I)\otimes V\otimes W \ar[dd]^{\id_{_{\hom(W,I)}}\otimes\, ev_{_{V,I}}\otimes\id_{_W}} \\ \\
\hom(V,I)\otimes V \ar[rr]_-{ev_{_{V,I}}} & & I & & \hom(W,I)\otimes W. \ar[ll]^-{ev_{_{W,I}}}
}\]
\end{lemma}
\noindent\emph{Proof:} Given an object $U$ in $\CC$, we will denote $U^*=\hom(U,I)$. In fact, we will prove that the following diagram is commutative, where the arrows that appear in the opposite direction are isomorphisms.
\[\tiny
\xymatrix{ V^*\otimes W^*\otimes W\otimes V
\ar[rr]^{ \id_{_{V^*\otimes W^*}} \otimes\, c_{_{_{W,V}}} } 
\ar[dd]^{ \id_{_{V^*}}\otimes\, ev_{_{_{W,I}}} \otimes\, \id_{_{_{V}}} }
\ar[ddrr]^{ \quad  \id_{_{V^*}} \otimes\, c_{_{W^*\otimes W, V}} }
& & 
V^*\otimes W^*\otimes V\otimes W
\ar[dd]^{ \id_{_{V^*}} \otimes\, c_{_{W^*,V}} \otimes\, \id_{_{_{W}}} }
& & 
W^*\otimes V^*\otimes V\otimes W
\ar[ll]_{ c_{_{W^*,V^*}} \otimes\, \id_{_{_{V\otimes W}}} } \ar[ddll]^{\quad \quad \quad c_{_{W^*,V^*\otimes V}} \otimes\, \id_{_{_{W}}} \quad } \ar[dd]_{ \id_{_{W^*}} \otimes\, ev_{_{_{V,I}}} \otimes\, \id_{_{_{W}}} } \\ \\
V^*\otimes I \otimes V \ar[dd]^{\simeq} \ar[dr]^{ \id_{{V^*}} \otimes \, c_{_{_{I,V}}} }
& &
V^*\otimes V\otimes W^*\otimes W \ar[dl]_{ \id_{_{V^*\otimes V}} \otimes\, ev_{_{_{W,I}}} \quad } \ar[dr]_{ ev_{_{_{V,I}}} \otimes\, \id_{_{W^*\otimes W}} \quad }
& &
W^*\otimes I\otimes W \ar[dd]^{\simeq} \ar[dl]_{ c_{_{W^*\otimes I}} \otimes\, \id_{_{_{W}}} \quad \quad } \\
&
V^*\otimes V\otimes I \ar[dl]^{\simeq}
\ar[r]_-{ ev_{_{_{V,I}}} \otimes\, \id_{_{_{I}}} }
& I\otimes I \ar[d]^{\simeq} &
I\otimes W^*\otimes W \ar[dr]^{\simeq}
\ar[l]^{ \id_{_{_{I}}} \otimes\, ev_{_{_{W,I}}} } &
\\
V^*\otimes V \ar[rr]_{ ev_{_{_{V,I}}} } & & I & & 
W^*\otimes W. \ar[ll]^{ ev_{_{_{W,I}}} }
}\]

Both triangles on top of the diagram commute by using the axioms in the definition of the braiding $c$, the quadrilaterals on both sides commute by naturality of $c$, and the remaining subdiagrams commute by the axioms of a monoidal category.
\qed

\begin{prop} \label{propbraidingxi}
Let $(\CC, \otimes, I, c)$ be a braided monoidal category, and let $V$ and $W$ be objects in $\CC$ such that $\hom(U,I)$ exists for $U=V,W,V\otimes W,W\otimes V$. The following diagram commutes:
\[
\xymatrix{
\hom(W,I)\otimes\hom(V,I) \ar[r]^-{\xi_{W,V}} \ar[d]_{c_{\hom(W,I),\hom(V,I)}} & \hom(V\otimes W,I) \ar[d]^{\hom\cc c_{_{W,V}},\id_{_I}\dd} \\
\hom(V,I)\otimes\hom(W,I) \ar[r]_-{\xi_{V,W}} & \hom(W\otimes V,I).}
\]
\end{prop}
\noindent\emph{Proof:}
\ Let us consider the following diagram in $\Vect$, which commutes by naturality of the adjunction isomorphism $\Theta$ of Definition \ref{defihomint}.
\[\footnotesize
\xymatrix{\Hom_\CC\cc \hom(W,I)\otimes\hom(V,I)\otimes V\otimes W,I \dd \ar[dddd]^{\Hom_\CC\cc \cc c_{_{\hom(W,I),\hom(V,I)}}\dd\inv \otimes c_{_{W,V}},\id_{_I} \dd} \ar[dr]^{\quad\quad\quad\quad\Theta_{\hom(W,I)\otimes\hom(V,I),V\otimes W,I}} & \\ & \Hom_\CC\cc \hom(W,I)\otimes\hom(V,I), \hom(V\otimes W,I) \dd \ar[dddd]_{\Hom_\CC\cc \cc c_{_{\hom(W,I),\hom(V,I)}}\dd\inv , \hom\cc c_{_{W,V}},\id_{_I}\dd \dd} \\ \\ \\
\Hom_\CC\cc \hom(V,I)\otimes\hom(W,I)\otimes W\otimes V,I \dd \ar[dr]_{\Theta_{\hom(V,I)\otimes\hom(W,I),W\otimes V,I \quad\quad\quad\quad}} & \\ &
\Hom_\CC\cc \hom(V,I)\otimes\hom(W,I), \hom(W\otimes V,I) \dd.
}\]
Using the previous lemma and chasing the element $ev_{_{W,I}}\circ \cc \id_{_{\hom(W,I)}}\otimes ev_{_{V,I}}\otimes\id_{_W} \dd$ in this diagram gives
\[\footnotesize
\xymatrix{ ev_{_{W,I}}\circ \cc \id_{_{\hom(W,I)}}\otimes ev_{_{V,I}}\otimes\id_{_W} \dd \ar@{|->}[ddd] \ar@{|->}[drrr] & & & \\ & & &
\xi_{W,V}  \ar@{|->}[ddd] \\ \\
ev_{_{V,I}}\circ \cc \id_{_{\hom(V,I)}}\otimes ev_{_{W,I}}\otimes\id_{_V} \dd \ar@{|->}[drrr] & & & \\ & & & \xi_{V,W},
}\]

from which precomposing with $c_{_{\hom(W,I),\hom(V,I)}}$ gives
\[\begin{array}{rcl}
\xi_{_{V,W}} \circ c_{_{\hom(W,I),\hom(V,I)}}
&=& \hom\cc c_{_{W,V}},\id_{_I}  \dd \circ \xi_{_{W,V}},
\end{array}
\]
which is what we wanted to prove. \qed

\medskip

\begin{rema} \label{remacheckcomm}
Let us consider a bialgebra $A$ in the braided monoidal category $\kGYD$ of Yetter-Drinfeld modules over a group algebra $\Bbbk G$. Suppose that either $A$ is finite dimensional or $G$ is finite, so that $\hom(A^{\otimes n}, \Bbbk) = \Homk (A^{\otimes n}, \Bbbk)$ - see Example 
\ref{exhomintkGYD}.
It follows from Theorem \ref{thmMPSW} that the diagram induced in cohomology by the following one commutes:
\[
\xymatrix{
\Homk(A^{\otimes p},\Bbbk)\otimes\Homk(A^{\otimes q},\Bbbk)  \ar[d]_{c_{_{\Homk\cc A^{\otimes p},\Bbbk  \dd,\Homk\cc A^{\otimes q},\Bbbk\dd}}} \ar[dr]^{(-1)^{pq}\smallsmile_{p,q}} & \\
\Homk(A^{\otimes q},\Bbbk)\otimes\Homk(A^{\otimes p},\Bbbk) \ar[r]_-{\smallsmile_{q,p}} & \Homk(A^{\otimes p+q},\Bbbk),}
\]
which, due to the previous proposition, is equivalent to the commutativity of the following diagram after taking cohomology:
\[
\xymatrix{
\Homk(A^{\otimes p},\Bbbk)\otimes\Homk(A^{\otimes q},\Bbbk)  \ar[r]^-{\smallsmile_{p,q}} \ar[dr]_{(-1)^{pq}\smallsmile_{p,q}} &
\Homk(A^{\otimes p+q},\Bbbk) \ar[d]^{\Homk\cc c_{A^{\otimes p},A^{\otimes q}},\Bbbk\dd} \\ & \Homk(A^{\otimes p+q},\Bbbk),}
\]

In order to check the commutativity of the first diagram, it is necessary to evaluate the braiding $c$ in the objects $\Homk(A^{\otimes p},\Bbbk)$ and $\Homk(A^{\otimes q},\Bbbk)$, which can be done because they coincide with their respective inner hom objects, so they can be regarded as objects in $\kGYD$.

On the other hand, checking the commutativity of the second diagram only requires evaluating the braiding in the objects $A^{\otimes p}$ and $A^{\otimes q}$. No inner hom objects are needed, so this condition can be checked even when the objects $\hom(A^{\otimes p},\Bbbk)$ and $\hom(A^{\otimes q},\Bbbk)$ do not exist, or when they exist but differ from $\Homk(A^{\otimes p},\Bbbk)$ and $\Homk(A^{\otimes p},\Bbbk)$, respectively. This is the case of the Jordan and the super Jordan plane. For these algebras, we have checked that the second diagram commutes at the level of cohomology, following the steps we give in detail in the Appendix.
\end{rema}

\medskip
In what follows, we will prove that the condition checked above holds up to homotopy for any bialgebra $A$ in a braided monoidal category $(\CC,\otimes,I,c)$, and give a categorical interpretation of the result.

For this purpose, we will perform some calculations using graphical calculus. Tensor products between objects are represented as their horizontal juxtaposition, and morphisms as sets of strings going from top to bottom. Here is a table with the notations we will use:
\medskip
\begin{center}
\begin{tabular}{| c | c | c | c |} 
 \hline
Morphism & Symbol & Morphism & Symbol\\  [1ex]
\hline 
&&&\\
 $c_{X,Y}: X\otimes Y \to Y\otimes X$ & \gbeg23 \got1X\got1Y\gnl \gbr\gnl \gob1Y\gob1X \gend & $\Delta: A \to A\otimes A$  & 
 
 $\gbeg23 \got2A\gnl \gcmu\gnl \gob1A\gob1A  \gend$ \\ &&& \\
 
 $\left(c_{Y,X}\right)^{-1}: X\otimes Y \to Y\otimes X$ & \gbeg23 \got1X\got1Y\gnl \gibr\gnl \gob1Y\gob1X \gend & 
 
 $\varepsilon: A\to I$ & $\gbeg13 \got1A\gnl \gcl1\gnl \gcu1 \gend$ \\ &&& \\

 $\mu: A\otimes A \to A$ & $\gbeg23 \got1A\got1A\gnl \gmu\gnl \gob2A \gend $ &
 
 $\chi_\ell: A\otimes M \to M$ & $\gbeg24
 \got1A\got1M\gnl
 \gcl1\gcl1\gnl 
 \glm\gnl 
 \gvac1\gob1M \gend$ \\ &&& \\

 $\eta : I\to A$ & $\gbeg13  \gu1\gnl \gcl1\gnl \gob1A  \gend$ &
 
 $\chi_r: M\otimes A \to M$ & $\gbeg24 \got1M\got1A\gnl \gcl1\gcl1\gnl \grm\gnl
 \gob1M \gend$ \\ &&& \\
 [1ex] 
 \hline
\end{tabular}
\end{center}
\medskip
These are examples of properties that we will use in the calculations: 
\begin{itemize}
\item Naturality of the braiding, as well as the axiom involving it and the unit constraints, can be expressed graphically when tensoring any morphism with a morphism from the unit object to another object of $\CC$ or vice versa:
\[
\gbeg24 \gvac1\got1Y\gnl \gu1\gcl1\gnl \gbr\gnl \gob1Y\gob1X\gend = 
\gbeg24 \got1Y\gnl \gcl1\gu1\gnl \gcl1\gcl1\gnl \gob1Y\gob1X\gend,\quad
\gbeg24 \got1X\gnl \gcl1\gu1\gnl \gbr\gnl \gob1Y\gob1X\gend =
\gbeg24 \gvac1\got1X\gnl \gu1\gcl1\gnl \gcl1\gcl1\gnl \gob1Y\gob1X\gend,\quad
\gbeg24 \got1X\got1Y\gnl \gbr\gnl \gcu1\gcl1\gnl \gvac1\gob1X\gend =
\gbeg24 \got1X\got1Y\gnl \gcl1\gcl1\gnl\gcl1\gcu1\gnl\gob1X\gend, \quad
\gbeg24 \got1X\got1Y\gnl \gbr\gnl \gcl1\gcu1\gnl \gob1Y\gend =
\gbeg24 \got1X\got1Y\gnl \gcl1\gcl1\gnl \gcu1\gcl1\gnl \gvac1\gob1Y\gend.
\]
\item One can also visualize naturality of the braiding in the following examples:
\[
\gbeg35 \got1A\got1A\got1X\gnl \gcl1\gbr\gnl \gbr\gcl1\gnl \gcl1\gmu\gnl \gob1X\gob2A\gend = 
\gbeg34 \got1A\got1A\got1X\gnl \gmu\gcl1\gnl \gbbrh3125 \gnl \gob2X\gob1A\gend, \quad\quad\quad
\gbeg35 \got2C\got1X\gnl \gcmu\gcl1\gnl \gcl1\gbr\gnl \gbr\gcl1\gnl \gob1X\gob1C\gob1C\gend =
\gbeg34 \got1C\got2X\gnl \gbbrh2114\gnl \gcl1\gcmu\gnl \gob1X\gob1C\gob1C\gend, \quad\quad\quad
\gbeg35 \got1A\got1M\got1X\gnl \gcl1\gbr\gnl \gbr\gcl1\gnl \gcl1\glm\gnl \gob1X\gob3M\gend = 
\gbeg34 \got1A\got1M\got1X\gnl \glm\gcl1\gnl \gbbrh3135 \gnl \gob2X\gob1M\gend.
\]

\item The following equalities represent the axioms of associativity and unit of a monoid:
\[
\gbeg54
\got1A \gvac1 \got1A \gvac1 \got1A \gnl
\gwmuh315 \gvac1 \gcl1 \gnl
\gwmuh639 \gnl
\gob6A
\gend = 
\gbeg54
\got1A \gvac1 \got1A \gvac1 \got1A \gnl
\gcl1 \gwmuh537 \gnl
\gwmuh417 \gnl
\gob4A
\gend,
\quad\quad
\gbeg24 \gvac1\got1A\gnl \gu1\gcl1\gnl \gmu\gnl \gob2A\gend =
\gbeg14 \got1A\gnl \gcl2\gnl\gnl \gob1A\gend =
\gbeg24 \got1A\gnl \gcl1\gu1\gnl \gmu\gnl \gob2A\gend,
\]
and the corresponding ones for a comonoid are the ones above turned upside down.
\end{itemize}

 \begin{ex}
 The simplicial object $(S(A),\{\partial_i\}_i,\{\sigma_j\}_j)$ introduced in Definition \ref{defiS} can be written using graphical notation as follows:
 \begin{itemize}
    \item Its objects are $S_n(A) = A\overset{n}{\cdots}A$,
    \item its faces $\partial_i\colon S_{n+1}(A)\to S_n(A)$ are given by
    
\[ \partial_0 = 
\gbeg45 \gvac2\got1{_n}\gnl \got1A\got1A\got1\cdots\got1A\gnl \gcl1\gcl1\gvac1\gcl1\gnl \gcu1\gcl1\got1\cdots\gcl1\gnl \gvac1\gob1A\gob1\cdots\gob1A\gend, \quad\quad
\partial_i = 
\gbeg65 \gvac1\got1{_{i-1}}\gvac2\got1{_{n-i}}\gnl \got1A\got1\cdots\got1A\got1A\got1\cdots\got1A\gnl \gcl1\gvac1\gcl1\gcl1\gvac1\gcl1\gnl \gcl1\got1\cdots\gmu\got1\cdots\gcl1\gnl \gob1A\gob1\cdots\gob2{\!\cdot   A\,\cdot}\gob1\cdots\gob1A\gend,\quad\quad
\partial_{n+1}
\gbeg45 \gvac1\got1{_n}\gnl \got1A\got1\cdots\got1A\got1A\gnl \gcl1\gvac1\gcl1\gcl1\gnl \gcl1\got1\cdots\gcl1\gcu1\gnl \gob1A\gob1\cdots\gob1A\gend,
\]
\item its degeneracies $\sigma_j\colon S_n(A)\to S_{n+1}(A)$ are given by
\[
\sigma_j = 
\gbeg75 \gvac1\got1{_{j-1}}\gvac3\got1{_{n+1-j}}\gnl \got1A\got1\cdots\got1A\gvac1\got1A\got1\cdots\got1A\gnl \gcl1\gvac1\gcl1\gu1\gcl1\gvac1\gcl1\gnl \gcl1\got1\cdots\gcl1\gcl1\gcl1\got1\cdots\gcl1\gnl \gob1A\gob1\cdots\gob1A\gob1A\gob1A\gob1\cdots\gob1A
\gend.
\]
\end{itemize}
\end{ex}
\medskip
\begin{defi}
Let $(S(A)\otimes S(A))_\bullet$ denote the product of $S_\bullet(A)$  with itself in $\Ch(\CC)$, defined analogously to Example \ref{defitensorCoch}. The \emph{deconcatenation coproduct} is the map $dec_\bullet: S_\bullet(A)\longrightarrow (S(A)\otimes S(A))_\bullet$ given componentwise by $dec_n=\sum_{p+q=n}dec_{p,q}$, where
\[
dec_{p,q}=\id_{A^{\otimes p+q}}: S_{p+q}(A) \longrightarrow S_p(A)\otimes S_q(A).
\]
\end{defi}

The following result is the generalization of the previous computations, which will allow us to recover the main result in Section 3 of \cite{mpsw}.
\begin{thm}\label{thmdec}
The triple $\left(S(A)_\bullet, dec_\bullet, (\id_{\Bbbk})_0\right)$ is a comonoid in $\Ch(\CC)$, which is cocommutative up to homotopy.
\end{thm}

In order to prove the theorem, we will introduce in our context the Alexander-Whitney map - defined more generally in \cite[Section 8.5]{weibel} - and recall Lemma 3.9 of \cite{mpsw}.

\begin{defi}
For an augmented monoid $A$ in a monoidal category $\CC$, let $S(A)\times S(A)$ denote the product of the simplicial object $S(A)$ in $\CC$ with itself, and $(S(A)\times S(A))_\bullet$ its associated chain complex. 
The \emph{Alexander-Whitney map} $\AW_\bullet: (S(A)\times S(A))_\bullet \longrightarrow (S(A)\otimes S(A))_\bullet$ is defined componentwise by the maps $\AW_{p,q}:S(A)_{p+q}\otimes S(A)_{p+q} \longrightarrow S(A)_p\otimes S(A)_q$,
\[
\AW_{p,q} = 
\gbeg{12}5 \gvac1\got1{_p}\gvac2\got1{_q}\gvac2\got1{_p}\gvac2\got1{_q}\gnl 
\got1A\got1\cdots\got1A\got1A\got1\cdots\got1A\got1A\got1\cdots\got1A\got1A\got1\cdots\got1A\gnl \gcl1\gvac1\gcl1\gcl1\gvac1\gcl1\gcl1\gvac1\gcl1\gcl1\gvac1\gcl1\gnl 
\gcl1\got1{\cdots}\gcl1\gcu1\got1{\cdots}\gcu1\gcu1\got1\cdots\gcu1\gcl1\got1\cdots\gcl1\gnl
\gob1A\gob1\cdots\gob1A\gvac6\gob1A\gob1{\cdots}\gob1A\gend.
\]

\medskip

The \emph{twisted Alexander-Whitney map} $\AWt_\bullet: (S(A)\times S(A))_\bullet \longrightarrow (S(A)\otimes S(A))_\bullet$ is defined componentwise by the maps $\AWt_{p,q}:S(A)_{p+q}\otimes S(A)_{p+q} \longrightarrow S(A)_p\otimes S(A)_q$,
\[
\AWt_{p,q} = (-1)^{pq}\
\gbeg{12}5 \gvac1\got1{_q}\gvac2\got1{_p}\gvac2\got1{_q}\gvac2\got1{_p}\gnl 
\got1A\got1\cdots\got1A\got1A\got1\cdots\got1A\got1A\got1\cdots\got1A\got1A\got1\cdots\got1A\gnl \gcl1\gvac1\gcl1\gcl1\gvac1\gcl1\gcl1\gvac1\gcl1\gcl1\gvac1\gcl1\gnl 
\gcu1\got1{\cdots}\gcu1\gcl1\got1{\cdots}\gcl1\gcl1\got1\cdots\gcl1\gcu1\got1\cdots\gcu1\gnl
\gvac3\gob1A\gob1\cdots\gob1A\gob1A\gob1{\cdots}\gob1A\gend.
\]
\end{defi}

As stated in \cite{weibel}, the maps $\AW_\bullet$ and $\AWt_\bullet$ are homotopic, since there is a unique natural chain map from $(S(A)\times S(A))_\bullet$ to $(S(A)\otimes S(A))_\bullet$ up to natural homotopy. Therefore, precomposing either of them with any chain map  $f_\bullet:S_\bullet(A) \to (S(A)\times S(A))_\bullet$ will give the same map up to homotopy. In the case $A$ is a bialgebra, we will use the map $f_\bullet$ as the composition of the functor $S_\bullet$ applied to the coproduct of $A$, and the map $g^{_{A,A}}_\bullet$ obtained from the simplicial map defined in the following lemma.

\begin{lemma} \cite[Lemma 3.9]{mpsw}
 Let $A$ and $B$ be two augmented algebras in a braided monoidal category $(\CC, \otimes, I, c)$. For each $n\in\N_0$, let us define recursively a map $g^{_{A,B}}_n: (A\otimes B)^{\otimes n} \longrightarrow A^{\otimes n}\otimes B^{\otimes n}$ as follows:
 \begin{itemize}
     \item $g^{_{A,B}}_0= \id_{_I}$,
     \item $g^{_{A,B}}_{n+1} = c_{_{B^{\otimes n},A}} \circ \cc g^{_{A,B}}_n\otimes \id_{_{A\otimes B}} \dd$
 \end{itemize}
The collection of maps $g^{_{A,B}}= \cc g^{_{A,B}}_n \dd_{n\in\N_0}$ is an isomorphism of simplicial objects in $\CC$ from $S(A\otimes B)$ to $S(A)\times S(B)$.
\end{lemma}

\medskip

\noindent {\bf Proof of Theorem \ref{thmdec}:}

Coassociativity and counitality of $\left(S(A)_\bullet, dec_\bullet, (\id_{\Bbbk})_0\right)$ are immediate consequences of the definition, and checking that the maps involved are chain morphisms is straightforward.

\noindent Proving cocommutativity up to homotopy is checking that the following diagram commutes up to homotopy:
\[
\xymatrix{ & S_\bullet(A) \ar[dl]_{dec\bullet} \ar[dr]^{dec_\bullet} & \\
(S(A)\otimes S(A))_\bullet \ar[rr]_{c^{gr}_{S_\bullet(A),S_\bullet(A)}} & & (S(A)\otimes S(A))_\bullet,
}
\]
which in the component $S_p(A)\otimes S_q(A)$ reads as follows:
\[
\xymatrix{ & S_{p+q}(A) \ar[dl]_{dec_{q,p}} \ar[dr]^{dec_{p,q}} & \\
S_q(A)\otimes S_p(A) \ar[rr]_{(-1)^{pq}\ c_{_{S_q(A),S_p(A)}}} & & S_p(A)\otimes S_q(A).
}
\]
Notice that
\[
\id_{A^{\otimes p+q}} = \AW_{p,q}\circ g^{_{A,A}}_{p+q}\circ S_{p+q}(\Delta) \sim \AWt_{p,q}\circ g^{_{A,A}}_{p+q}\circ S_{p+q}(\Delta) = (-1)^{pq}c_{A^{\otimes q}, A^{\otimes p}},
\]
where the equalities are due to counitality. To illustrate them, we show an example for $p=2,q=1$. In the second and third diagrams, the level at the top corresponds to the morphism $S_3(\Delta)$, the middle ones all together correspond to the morphism $g_3^{A,A}$, and the one at the bottom corresponds to the morphisms $\AW_{2,1}$ and $\AWt_{2,1}$, respectively.
\[
\gbeg36
    \got1A\got1A\got1A\gnl
    \gcl4\gcl4\gcl4\gnl\gnl
    \gob1A\gob1A\gob1A
\gend
\quad = \quad
\gbeg66
    \got2A\got2A\got2A\gnl
    \gcmu\gcmu\gcmu\gnl
    \gcl2\gbr\gbr\gcl2\gnl
    \gvac1\gcl1\gbr\gcl1\gnl
    \gcl1\gcl1\gcu1\gcu1\gcu1\gcl1\gnl
    \gob1A\gob1A\gob7A
\gend
\quad \sim \quad
(-1)^{2\cdot 1}\
\gbeg66
    \got2A\got2A\got2A\gnl
    \gcmu\gcmu\gcmu\gnl
    \gcl2\gbr\gbr\gcl2\gnl
    \gvac1\gcl1\gbr\gcl1\gnl
    \gcu1\gcl1\gcl1\gcl1\gcu1\gcu1\gnl
    \gvac1\gob1A\gob1A\gob1A
\gend
\quad = \quad
(-1)^{2\cdot 1}\
\gbeg36
    \got1A\got1A\got1A\gnl
    \gcl1\gcl1\gcl1\gnl
    \gbr\gcl1\gnl
    \gcl1\gbr\gnl
    \gcl1\gcl1\gcl1\gnl
    \gob1A\gob1A\gob1A
\gend 
\]
\qed

To conclude this section, we will interpret Proposition \ref{propbraidingxi} in a categorical way, using the definition of braided monoidal functors. This will allow us to recover the result of Corollary \ref{corMPSW}.

\begin{defi}
Let $(\CC,\otimes,I,c^\CC)$ and $(\DD,\times,J,c^\DD)$ be two braided monoidal categories. A lax monoidal functor  $(F,\varphi,\varphi_0):(\CC,\otimes,I) \to (\DD,\times,J)$ is a \emph{braided lax monoidal functor} if for every $X,Y\in\CC$ the following diagram commutes:
\[
\xymatrix{ F(X)\times F(Y) \ar[rr]^{c^\DD_{F(X),F(Y)}} \ar[d]_{\varphi_{X,Y}} & & F(Y)\times F(X) \ar[d]^{\varphi_{Y,X}} \\ F(X\otimes Y) \ar[rr]_{F\left(c^\CC_{X,Y}\right)} & & F(Y\otimes X).
}
\]
\end{defi}

\begin{ex}
The braiding defined in Example \ref{extau} on the categories $^H\Mod$ for a commutative $\Bbbk$-bialgebra $H$ and $_H\Mod$ for a cocommutative $\Bbbk$-bialgebra $H$ coincides with the one of their underlying vector spaces, which automatically implies that their respective forgetful functors to the category of $\Bbbk$-vector spaces are braided.
\end{ex}

\begin{ex}
Unlike the previous examples, the forgetful functors corresponding to the categories $\Vect^\Z,\Ch(\Vect),\Coch(\Vect)$, and $\HYD$ for a Hopf algebra $H$ are not braided.
\end{ex}

\begin{ex}\label{exxibraidedmon}
This example provides a tool to recover Corollary \ref{corMPSW} from Theorem \ref{thmdec}.

Let $(\CC, \otimes, I)$ be a monoidal category, and let $Adj(\CC)$ be the full subcategory of $\CC$ whose objects are the $V\in\CC$ such that the functor $-\otimes V$ has right adjoint $\hom(V,-)$.

If $V$ and $W$ are objects in $Adj(\CC)$, the composition of natural isomorphisms
\[
\Hom_\CC(-\otimes V\otimes W, -) \simeq \Hom_\CC(-\otimes V, \hom(W,-)) \simeq \Hom_\CC(-, \hom(V,\hom(W,-)))
\]
implies that $V\otimes W$ is an object of $Adj(\CC)$, and $\hom(V\otimes W, Y) = \hom(V, \hom(W,Y))$ for every object
$Y$ of $\CC$.
In addition, the natural isomorphism
\[
\Hom_\CC(- \otimes I, -) \simeq \Hom_\CC(-,-)
\]
implies that $I\in Adj(\CC)$ and $\hom(I, Y) = Y$ for every object $Y\in\CC$. Hence, one can consider the monoidal category $(Adj(\CC), \otimes, I)$ and the restriction of the functor $\hom(-,I):Adj(\CC)^\flop \to\CC$, which has a monoidal structure 
\[
(\hom(-,I), \xi, \xi_0): (Adj(\CC)^\flop, \totimes, I) \to (\CC, \otimes, I),
\]
Here $\xi$ is the natural transformation of Definition \ref{defixi}, and $\xi_0=\id_I$. Proposition \ref{propbraidingxi} is the condition for this lax monoidal functor to be braided.
\end{ex}

The last step is applying this proposition to the lax monoidal functor defined above.
\begin{prop} \label{propbraidedsendscomm} \cite[Proposition 3.37]{aguiar}
Given two braided monoidal categories $(\CC,\otimes,I)$ and $(\DD,\times,J)$, let $(F,\varphi,\varphi_0)$ be a braided lax monoidal functor from $\CC$ to $\DD$, and let $(A,\mu,\eta)$ be a commutative monoid in $\CC$. The triple $\left(F(A),F(\mu)\circ\varphi_{A,A},F(\eta)\circ\varphi_0\right)$ is a commutative monoid in $\DD$.
\end{prop}

We finish this section with a result of graded braided commutativity up to homotopy in $\hCoch\left(\HYD\right)$. It slightly generalizes Corollary 3.13 of \cite{mpsw}, which can be obtained from this theorem by taking cohomology.

\begin{thm} \label{thmMPSWhomotopy}
Let $H$ be a Hopf algebra, and let $A$ be a bimonoid in the braided monoidal category $\left(\HYD, \otimes, \Bbbk, c\right)$. Suppose that either $H$ or $A$ is finite dimensional. Then the differential graded algebra $\Homk(S(A)_\bullet,\Bbbk)$ with the opposite of the cup product is graded braided commutative up to homotopy.
\end{thm}
\noindent\emph{Proof:} \ Let us consider the chain complex $S_\bullet(A)$ with the deconcatenation coproduct $dec$ and its counit. By Theorem \ref{thmdec}, it is a graded braided cocommutative comonoid in $\hCh\left(\HYD\right)$.

On the other hand, the finite dimension of either $H$ or $A$ guarantees that the objects $S_n(A)$ lie in the subcategory
$Adj\left(\HYD\right)$ defined in Example \ref{exxibraidedmon}. Let us denote this subcategory by $\BB$. One can think of $S_\bullet(A)$ as a chain complex in $\BB$.
The functor $\hom(-,I)$ applied to the objects and the differentials of a chain complex, allows us to define a functor $\hom(-,I)^\Coch: \cc\Ch(\BB)\dd^\flop \simeq \Coch(\BB^\flop) \longrightarrow \Coch(\CC)$.

One can now define componentwise a lax monoidal functor 
\[
\left(\hom(-,I)^\Coch,\xi^\Coch, \xi^\Coch_0\right): \left(\Coch(\BB^\flop),\Tilde{\otimes},I_0\right) \longrightarrow \left(\Coch(\CC),\otimes,I_0\right),
\]
where the monoidal structures in the categories of cochain complexes are those of Remark \ref{remaChCochCC}, making $\xi^\Coch_{V_\bullet,W_\bullet}$ a morphism of cochain complexes. Associativity and unitality are checked componentwise, as well as the fact that this lax monoidal functor is braided with the graded braiding on both categories.

The same construction can be made modulo homotopy. Since $(S_\bullet(A),dec,\id_\Bbbk)$ is a commutative monoid in the braided monoidal category $(\hCoch(\BB^\flop),\Tilde{\otimes},\Bbbk_0, \Tilde{c}^{gr})$, Proposition \ref{propbraidedsendscomm} gives a commutative monoid in $(\hCoch(\CC),\otimes,\Bbbk_0, c^{gr})$, which coincides as a cochain complex with $\hom(S_\bullet(A),\Bbbk)$. It is straightforward to check that the product of Proposition \ref{propbraidedsendscomm} coincides with the cup product $\smallsmile$. 

Finiteness of the dimension of either $H$ or $A$ also guarantees that the vector spaces $\hom(S_n(A),\Bbbk)$ and $\Homk(S_n(A),\Bbbk))$ are isomorphic, and that the product $\smallsmile$ is the opposite of the usual cup product $\cup$.

Therefore, $\Homk(S_\bullet(A),\Bbbk)$ with the opposite of the usual cup product is a graded braided commutative algebra up to homotopy, and so it is after taking cohomology.
\qed

\section{An approach using duoidal categories}

In the previous section we proved a cocommutativity theorem which holds in several braided monoidal categories, and which is also independent of the existence of cochain complexes and cohomology inside the category. However, it is still desirable to have commutativity at the level of cohomology. For this reason, our aim in this section is to show a result of this kind for a more general class of algebras than in \cite{mpsw}.

One problem of the previous approach is that it depends on the bar resolution, which is in general difficult to handle. Although some computations can be made, it would be better to be able to work with smaller projective resolutions. There is a technical issue when defining the comultiplication of a projective resolution: it is well-defined only on the tensor product over the algebra $A$, which in general admits no braidings. Instead, it is possible to deal with it using a generalization of braided monoidal categories: duoidal categories.

\begin{defi} \cite[Definition 6.1: \emph{2-monoidal category}]{aguiar} A \emph{duoidal category} is a 15-tuple $\left(\CC, \diamond, I, a, \ell, r, \star, J, \alpha, \lambda, \rho, \zeta, \Delta_I, \mu_J, \zeta_0\right)$, where 
\begin{itemize}
\item $(\CC, \diamond, I, a, \ell, r)$ and $(\CC, \star, J, \alpha, \lambda, \rho)$ are monoidal categories.
\item $(I, \Delta_I, \zeta_0)$ is a comonoid in $(\CC, \star, J, \alpha, \lambda, \rho)$.
\item $(J, \mu_J, \zeta_0)$ is a monoid in $(\CC, \diamond, I, a, \ell, r)$
\item $\zeta$ is a natural transformation, called \emph{interchange law}, given by a family of morphisms $\zeta_{X,Y,Z,T}: (X\star Y) \diamond (Z\star T) \to (X\diamond Z) \star (Y\diamond T)$ for each $X,Y,Z,T\in\CC$, such that the following diagrams commute for every $X,Y,Z,T,U,V \in\CC$:
\end{itemize}
\[
\xymatrix{
((X\star T)\diamond (Y\star U)) \diamond (Z\star V)
\ar[rr]^{a_{_{X\star T, Y\star U, Z\star V}}}
\ar[d]_{\zeta_{_{X,T,Y,U}}\diamond\,\id_{_{Z\star V}}} & &
(X\star T)\diamond ((Y\star U) \diamond (Z\star V)) 
\ar[d]^{\id_{_{X\star T}}\diamond\,\zeta_{_{Y,U,Z,V}}} \\
((X\diamond Y)\star (T\diamond U)) \diamond (Z\star V)
\ar[d]_{\zeta_{_{X\diamond Y, T\diamond U, Z,V}}} & &
(X\star T) \diamond ((Y\diamond Z)\star(U\diamond V)) 
\ar[d]^{\zeta_{_{X,T, Y\diamond Z, U\diamond V}}} \\
((X\diamond Y)\diamond Z)\star ((T\diamond U)\diamond V)
\ar[rr]_{a_{_{X,Y,Z}}\diamond\,a_{_{T,U,V}}} & &
(X\diamond (Y \diamond Z))\star (T\diamond (U\diamond V)), \\}
\]
\medskip
\[
\xymatrix{
((X\star Y)\star Z) \diamond ((T\star U)\star V) 
 \ar[rr]^{\alpha_{_{X,Y,Z}}\diamond\,\alpha_{_{T,U,V}}} 
 \ar[d]_{\zeta_{_{X\star Y, Z, T\star U,V}}} & &
(X\star (Y \star Z)) \diamond (T\star (U\star V))
\ar[d]^{\zeta_{_{X, Y\star Z, T, U\star V}}} \\
((X\star Y)\diamond (T\star U)) \star (Z\diamond V) 
\ar[d]_{\zeta_{_{X,T,Y,U}}\star\,\id_{_{Z\diamond V}}} & & (X\diamond T) \star ((Y\star Z)\diamond(U\star V))
\ar[d]^{\id_{_{X\diamond T}}\star\,\zeta_{_{Y,Z,U,V}}} \\
((X\diamond T)\star (Y\diamond U)) \star (Z\diamond V)
\ar[rr]_{\alpha_{_{X\diamond T, Y\diamond U, Z\diamond W}}} & &
(X\diamond T)\star ((Y\diamond U) \star (Z\diamond V)),}
\]

\[
\xymatrix{
I\diamond (X\star Y) \ar[r]^{\ell_{_{X\star Y}}} \ar[d]_{\Delta_{_{I}}\diamond\,\id_{_{X\star Y}}}  & X\star Y \\
(I\star I)\diamond (X\star Y) \ar[r]_{\zeta_{_{I,I,X,Y}}} & (I\diamond X)\star (I\diamond Y) \ar[u]_{\ell_{_{X}}\star\,\ell_{_{Y}}} ,}
\qquad
\xymatrix{
(X\star Y) \diamond I \ar[r]^{r_{_{X\star Y}}} \ar[d]_{\id_{_{X\star Y}}\diamond\, \Delta_{_{I}}} & X\star Y \\
(X\star Y)\diamond (I\star I) \ar[r]_{\zeta_{_{X,Y,I,I}}} & (X\diamond I)\star (Y\diamond I) \ar[u]_{r_{_{X}}\star\, r_{_{Y}}} ,}
\]
\[
\xymatrix{
(J\star X)\diamond (J\star Y) \ar[r]^{\zeta_{_{J,X,J,Y}}} \ar[d]_{\lambda_{_{X}}\diamond\,\lambda_{_{Y}}} &
(J\diamond J)\star (X\diamond Y) \ar[d]^{\mu_{_{J}}\star\, \id_{_{X\diamond Y}}}
\\ X\diamond Y & J\star (X\diamond Y) \ar[l]^{\lambda_{_{X\diamond Y}}} ,}
\qquad
\xymatrix{
(X\star J)\diamond (Y\star J) \ar[r]^{\zeta_{_{X,J,X,J}}} \ar[d]_{\rho_{_{X}}\diamond\,\rho{_{Y}}} &
(X\diamond Y)\star (J\diamond J) \ar[d]^{\id_{_{X\diamond Y}}\star\,\mu_{_{J}}}
\\ X\diamond Y &  (X\diamond Y)\star J \ar[l]^{\rho{_{X\diamond Y}}} ,}
\]

We will omit in general the structure on the units, denoting the duoidal category by $(\CC,\diamond, I, \star, J, \zeta, \zeta_0)$.
\end{defi}

\begin{ex}
    Any braided monoidal category $(\CC,\otimes,I,a,\ell,r,c)$ gives rise to the duoidal category
 $\cc\CC, \otimes, I, a, \ell, r, \otimes, I, a, \ell, r, \zeta^c, \ell_I\inv=r_I\inv,\ell_I=r_I,\id_I\dd$, where the interchange law is $\zeta^c_{X,Y,Z,T}=\id_X\otimes c_{Y,Z}\otimes\id_T$.
\end{ex}
\begin{rema}
Note that all morphisms involved in the previous example are isomorphisms. There is a categorical Eckmann-Hilton argument which states that essentially every duoidal category whose structure morphisms are isomorphisms comes from a braided monoidal category. See, for example, \cite[Proposition 6.11]{aguiar}.
\end{rema}

\begin{rema}
If $\left(\CC, \diamond, I, a, \ell, r, \star, J, \alpha, \lambda, \rho, \zeta, \Delta_I, \mu_J, \zeta_0\right)$ is a duoidal category, then so is its opposite category
 $\left(\CC^\flop, \star, J, \alpha\inv, \lambda\inv, \rho\inv, \diamond, I, a\inv, \ell\inv, r\inv, \zeta^\op, \mu_J, \Delta_I, \zeta_0\right)$, where the interchange law is $\zeta^\op_{X,Y,Z,T}=\zeta_{X,Z,Y,T}$.
\end{rema}

\begin{rema}
Transposing one or both of the monoidal structures of a monoidal category If $\left(\CC, \diamond, I, a, \ell, r, \star, J, \alpha, \lambda, \rho, \zeta, \Delta_I, \mu_J, \zeta_0\right)$ gives three different duoidal categories:
\begin{itemize}
    \item $\left(\CC, \tdiamond, I, a\inv, r, \ell, \star, J, \alpha, \lambda, \rho, \zeta^\tdiamond, \Delta_I, \mu_J, \zeta_0\right)$, where $\zeta^\tdiamond_{X,Y,Z,T}=\zeta_{Z,T,X,Y}$,
    \item $\left(\CC, \diamond, I, a, \ell, r, \tstar, J, \alpha\inv, \rho, \lambda, \zeta^\tstar, \Delta_I, \mu_J, \zeta_0\right)$, where $\zeta^\tstar_{X,Y,Z,T}=\zeta_{Y,X,T,Z},$ 
    \item $\left(\CC, \tdiamond, I, a\inv, r, \ell, \tstar, J, \alpha\inv, \rho, \lambda, \zeta^{\tdiamond,\tstar}, \Delta_I, \mu_J, \zeta_0\right)$, where $\zeta^{\tdiamond,\tstar}_{X,Y,Z,T}=\zeta_{T,Z,Y,X}$.
\end{itemize}
\end{rema}

\begin{defi}
A duoidal category is \emph{strict} if both underlying monoidal structures are strict.
\end{defi}

All duoidal categories considered in this work will be strict.

\begin{rema}\label{remazetagr}
For a duoidal category $\CC$, both monoidal structures can be lifted to $\Ch(\CC)$, $\Coch(\CC)$ and $\CC^\Z$, with the structures on the units concentrated in degree $0$ and the interchange law $\zeta^{gr}$ defined by $\zeta^{gr}_{X,Y,Z,T}|_{(X_m\star Y_n)\diamond (Z_k\star T_\ell)} = (-1)^{nk}\ \zeta_{X_m,Y_n,Z_k,T_\ell}$.
\end{rema}
\medskip

Next we will exhibit the duoidal category we are interested in.

\begin{defi}\label{defibimod}

\begin{enumerate}
\item Let $(\CC,\otimes,I)$ be a monoidal category, and let $(A,\mu,\eta)$ be a monoid in $\CC$. An \emph{$A$-bimodule} is a triple $(M,\chi_\ell, \chi_r)$, where
\begin{itemize}
    \item $M$ is an object in $\CC$,
    \item $\chi_\ell:A\otimes M\to M$ and $\chi_r:M\otimes A\to M$ are morphisms in $\CC$, respectively called \emph{left} and \emph{right action}, such that the following diagrams commute:
    \end{itemize}
    \[\xymatrix{
    A\otimes A\otimes M \ar[rr]^{\mu\otimes\id_M} \ar[d]_{\id_A\otimes\chi_\ell}& & A\otimes M \ar[d]^{\chi_\ell} \\
    A\otimes M \ar[rr]_{\chi_\ell} & & M,
    } \qquad
    \xymatrix{
    M\otimes A\otimes A \ar[rr]^{\chi_r\otimes\id_A} \ar[d]_{\id_M\otimes\mu}& & M\otimes A \ar[d]^{\chi_r} \\
    M\otimes A \ar[rr]_{\chi_r} & & M,
    }\]
    
    \[\xymatrix{
    A\otimes M\otimes A \ar[rr]^{\chi_\ell\otimes\id_A} \ar[d]_{\id_A\otimes\chi_r}& & M\otimes A \ar[d]^{\chi_r} \\
    A\otimes M \ar[rr]_{\chi_\ell} & & M,
    }\qquad\xymatrix{
    I\otimes M \ar[r]^{\eta\otimes\id_M} \ar[dr]_{\id_M} & A\otimes M \ar[d]^{\chi_\ell} \\
     & M.
    }\qquad
    \xymatrix{
    M\otimes A \ar[d]_{\chi_r} & M\otimes I \ar[l]_{\id_M\otimes\eta} \ar[dl]^{\id_M} \\
     M. &
    }\]

\item Let $(\CC,\otimes,I)$ be a monoidal category, and let $(A,\mu,\eta)$ be a monoid in $\CC$. The \emph{category $\ACA$ of $A$-bimodules in $\CC$} is defined as follows:
\begin{itemize}
    \item Objects are $A$-bimodules $(M, \chi_\ell,\chi_r)$.
    \item For each pair of $A$-bimodules $(M, \chi_\ell,\chi_r)$ and $(M', \chi'_\ell,\chi'_r)$, the set of morphisms $\Hom_\ACA(M,M')$ consists of the maps $f:M\to M'$ in $\CC$ such that the following diagrams commute:
\[
\xymatrix{
A\otimes M \ar[r]^-{\chi_\ell} \ar[d]_{\id_A\otimes f} & M \ar[d]^{f} \\
A\otimes M' \ar[r]_-{\chi'_\ell}& M',
}
\qquad
\xymatrix{
M\otimes A \ar[r]^-{\chi_r} \ar[d]_{f\otimes\id_A} & M \ar[d]^{f} \\
M'\otimes A \ar[r]_-{\chi'_r}& M'.
}
\]
\end{itemize}
\end{enumerate}
\end{defi}

Although mentioning $\CC$ is redundant in this definition, it will be necessary when dealing at the same time with different categories with the same underlying monoidal structures (e.g. vector spaces, graded vector spaces, Yetter Drinfeld modules).

\begin{defi}
Under the conditions of Definition \ref{defibimod}, let $(M, \chi_\ell,\chi_r)$ and $(M', \chi'_\ell,\chi'_r)$ be two $A$-bimodules in $\CC$. The \emph{product over $A$ of $M$ and $M'$} is the coequalizer $\pi_{M,M'}: M\otimes M' \to M\otimes_A M'$ of the maps  $\chi_r\otimes_{M'}$ and $\id_M\otimes\chi'_\ell:M\otimes A\otimes M'\longrightarrow M\otimes M'$. The left and right actions of $A$ on $M\otimes_A M'$ are respectively induced by $\chi_\ell$ and $\chi'_r$.
\end{defi}

\begin{rema}
It is easy to see that the product $\otimes_A$ is associative with unit $A$ and functorial, therefore $(\ACA, \otimes_A, A)$ is a monoidal category.
\end{rema}

\begin{defi}\label{defiodot}
Let $(\CC,\otimes,I,c)$ be a braided monoidal category, let $(A,\mu,\eta,\Delta,\varepsilon)$ be a bimonoid in $\CC$, and let $(M, \chi_\ell,\chi_r)$ and $(M', \chi'_\ell,\chi'_r)$ be $A$-bimodules in $\CC$. The $A$-bimodule  $M\odot M'$ is $\cc M\otimes M', \chi^\odot_\ell, \chi^\odot_r\dd$, with actions are given by
\[
\begin{array}{rcl}
\chi^\odot_\ell &=& \cc \chi_\ell \otimes \chi'_\ell \dd \circ \cc \id_A \otimes c_{A,M} \otimes \id_{M'} \dd \circ \cc \Delta \otimes \id_M \otimes \id_{M'} \dd, \\
\chi^\odot_r &=& \cc \chi_r \otimes \chi'_r \dd \circ \cc \id_M \otimes c_{M',A} \otimes \id_A \dd \circ \cc \id_M \otimes \id_{M'} \otimes \Delta \dd.
\end{array}
\]
\end{defi}
\begin{rema}
It is clear that the product of two $A$-bimodule morphisms is also an $A$-bimodule morphism, so $\odot$ is a bifunctor. Coassociativity and counitality of $A$ and the axioms of the braiding imply that $\odot$ is associative with unit $(I,\varepsilon,\varepsilon)$. Therefore, $(\ACA,\odot,I)$ is a monoidal category.
\end{rema}

We now obtain a duoidal category that will be used in what follows.

\begin{prop}
Let $(\CC,\otimes,I,c)$ be a braided monoidal category, let $(A,\mu,\eta,\Delta,\varepsilon)$ be a bimonoid in $\CC$. The 9-tuple $\cc \ACA,\otimes_A,A,\odot,I,\zeta^A , \Delta, \id_I, \varepsilon\dd$ is a duoidal category, where 
\[
\zeta^A_{M,N,K,L}: (M\odot N) \otimes_A (K\odot L) \longrightarrow (M\otimes_A K) \odot (N\otimes_A L)
\]
is the canonical projection of the map $\zeta^c_{M,N,K,L} = \id_M\otimes c_{N,K}\otimes \id_L$.
\end{prop}
\noindent \emph{Proof:} Compatibility of $\zeta^A$ with associativity and units of both products is a consequence of the axioms of the braiding $c$. We only need to prove that $\zeta^A$ is well-defined, that is, the map $\cc \pi_{M,K}\otimes \pi_{N,L} \dd \circ \zeta^c_{M,N,K,L}$ factors through $(M\odot N) \otimes_A (K\odot L)$. This is shown in the following graphical calculation. On one side, precomposing $\zeta^c_{M,N,K,L}$ with the right action of $A$ on $M\odot N$ yields:
\[
\gbeg66 \got1M \got1N \got2A \got1K \got1L \gnl
\gcl2 \gcl1 \gcmu \gcl3 \gcl4 \gnl
\gvac1 \gbr \gcl1 \gnl
\grm \grm \gnl\gcl1 \gvac1 \gbbrh2115 \gnl
\gob1M \gob3K \gob1N \gob1L \gend
=
\gbeg66 \got1M \got1N \got2A \got1K \got1L \gnl
\gcl3 \gcl1 \gcmu \gcl1 \gcl4 \gnl
\gvac1 \gbr \gbr \gnl
\gvac1 \gcl1 \gbr \gcl1 \gnl
\grm \gcl1 \grm \gnl 
\gob1M \gob3K \gob3{L.} \gnl
\gvac3 \gob1N \gend
\]
On the other side, precomposing $\zeta^c_{M,N,K,L}$ with the left action of $A$ on $P\odot L$ yields:
\[
\gbeg66 \got1M \got1N \got2A \got1K \got1L \gnl
\gcl4 \gcl3 \gcmu \gcl1 \gcl2 \gnl
\gvac2 \gcl1 \gbr  \gnl
\gvac2 \glm \glm \gnl
\gvac1 \gbbrh2115 \gvac2 \gcl1 \gnl
\gob1M \gob1K \gob3N \gob1L \gend
=
\gbeg66 \got1M \got1N \got2A \got1K \got1L \gnl
\gcl4 \gcl1 \gcmu \gcl1 \gcl3 \gnl
\gvac1 \gbr \gbr \gnl
\gvac1 \gcl1 \gbr \gcl1 \gnl
\gvac1 \glm \gcl1 \glm \gnl 
\gob1M \gob3K \gob3{L} \gnl
\gvac3 \gob1N \gob3{\quad ,} \gend
\]
which equals the previous calculation when postcomposed with $\pi_{M,K}\otimes \pi_{N,L}$. \qed

\medskip

This result holds in general for a duoidal category $\CC$ not necessarily arising from a braided monoidal category - see \cite[Section 7.1]{garnerlopez}. 

\medskip

The following definition generalizes the notion of a commutative monoid in a braided monoidal category to the context of duoidal categories.

\begin{defi}
Let $(\CC,\diamond, I, \star, J, \zeta, \Delta_I, \mu_J, \zeta_0)$ be a duoidal category. A \emph{duoid} in $\CC$ is a 5-tuple $(A,\mu,\eta,\nu,\iota)$, where
\begin{itemize}
    \item $(A,\mu,\eta)$ is a monoid in $(\CC,\diamond, I)$
    \item $(A,\nu, \iota)$ is a monoid in $(\CC,\star, I)$
    \item the following diagrams commute:
\[
\xymatrix{
 (A\star A)\diamond (A\star A) \ar[rr]^{\zeta_{A,A,A,A}} \ar[d]_{\nu\diamond\nu} & & (A\diamond A)\star (A\diamond A) \ar[d]^{\mu\star\mu} \\
 A\diamond A \ar[r]_-\mu & A & A\star A, \ar[l]^-\nu
}\]
\[
\xymatrix{
I \ar[r]^-{\Delta_I} \ar[d]_{\eta} & I\star  I \ar[d]^{\eta\star\eta} \\
A & A\star A \ar[l]^{\nu},
} \quad 
\xymatrix{
J\diamond J \ar[r]^-{\mu_J} \ar[d]_{\iota\diamond\iota}  & J \ar[d]^\iota \\
A\diamond A \ar[r]_-{\mu}& A,
} \quad
\xymatrix{
I \ar[rr]^{\zeta_0} \ar[dr]_\eta & & J \ar[dl]^\iota\\ & A. &
}
\]
\end{itemize}
\end{defi}

Also, for a duoidal category $(\CC,\diamond,I,\star,J,\zeta)$, a \emph{coduoid} in $\CC$ is a duoid in its opposite category $\CC^\flop$.

\begin{ex}\label{exEckmannHilton}
A commutative monoid $(A,\mu,\eta)$ in a braided monoidal category $(\CC,\otimes,I,c)$ gives a duoid $(A,\mu,\eta,\mu,\eta)$ in the duoidal category $(\CC,\otimes,I,\otimes,I,\zeta^c)$, and analogously a cocommutative comonoid induces a duoid. The converse is also true, due to the Eckmann-Hilton argument - see \cite[Proposition 6.29]{aguiar}.
\end{ex}
\medskip

In the sequel, we will focus on giving a coduoid structure to certain projective resolutions of $\Bbbk$ as $A$-bimodule, when $A$ is a bimonoid in the category $\kGYD$ of Yetter-Drinfeld modules over an abelian group $G$. One of the coproducts will be the one of Example \ref{exomega}, which is related to the definition of the usual cup product - see Remark \ref{remamultbimod}. First we will define the duoidal category whose coduoids we are interested in.

\begin{defi}
Let $A$ be a bimonoid in the category $\kGYD$. We define the duoidal category $\left(\khCh\left(\AkGYDA\right), \otimes_A, A_0, \odot, \Bbbk_0, \left( \zeta^A \right)^{gr} \right)$ as follows:
\begin{itemize}
    \item Objects are chain complexes of $A$-bimodules in $\kGYD$.
    \item Morphisms are $\Bbbk$-linear homotopy classes of chain morphisms in $\AModA$.
    \item The products $\otimes_A$ and $\odot$ are defined componentwise as in Remark \ref{remaChCochCC} with their units concentrated in degree $0$.
    \item The graded interchange law $\left( \zeta^A \right)^{gr}$ is defined as in Remark \ref{remazetagr}.
\end{itemize}
\end{defi}

\begin{thm} \label{thmPcoduoid}
Let $A$ be a bimonoid in the braided monoidal category $\kGYD$ for a group $G$, and let $P_\bullet \xrightarrow{f} A$ be a chain complex in $\AkGYDA$ which is also a projective resolution of $A$ in $\AModA$. The map $\omega:P_\bullet \to (P\otimes_A P)_\bullet$ lifting the isomorphism $A\simeq A\otimes_A A$ as in Example \ref{exomega} and the map $\delta: P_\bullet \to (P\odot P)_\bullet$ defined as a lifting of the map $\Delta:A \to A\odot A$ provide a coduoid $(P_\bullet, \omega, f_0, \delta, \varepsilon \circ f_0)$ in $\left(\khCh\left(\AkGYDA\right), \otimes_A, A_0, \odot, \Bbbk_0, \left( \zeta^A \right)^{gr} \right)$.
\end{thm}
\noindent\emph{Proof:} Consider the following commutative diagram in $\AkGYDA$:
\[
\xymatrix{ A\otimes_A A \ar[d]_{\Delta \otimes_A \Delta} & A \ar[l]_-\simeq \ar[r]^-\Delta & A\odot A \ar[d]^\simeq  \\
 (A\odot A)\otimes_A (A\odot A) \ar[rr]_{\zeta^A_{A,A,A,A}}  & & (A\otimes_A A)\odot (A\otimes_A A).
 }
\]
Since the forgetful functor $\kGYD\to\Vect$ is strong monoidal, the diagram is also commutative in $\AModA$. Lifting these morphisms along the exact complexes in $\AModA$ involving products of $P_\bullet$ gives the following diagram in $\AModA$:
\[
\xymatrix{ (P\otimes_A P)_\bullet \ar[d]_{\delta \otimes_A \delta} & P_\bullet \ar[l]_-\omega \ar[r]^-\delta & (P\odot P)_\bullet \ar[d]^{\omega \odot \omega}  \\
 ((P\odot P)\otimes_A (P\odot P))_\bullet \ar[rr]_{\left(\zeta^A\right)^{gr}_{P_\bullet,P_\bullet,P_\bullet,P_\bullet}}  & & ((P\otimes_A P)\odot (P\otimes_A P))_\bullet,
 }
\]
which consequently commutes up to a $\Bbbk$-linear homotopy.
\qed

\medskip

In the next examples we show how the Jordan and super Jordan plane fit into the hypotheses of the previous theorem.

\begin{ex} \label{exresoljordan}
For the Jordan plane, the following resolution of $A$ as $A$-bimodule has been computed in \cite{lopessolotar}:

\[\xymatrix{
0 \ar[r] & A\otimes R\otimes A \ar[r]^-{d_1} & A\otimes V\otimes A \ar[r]^-{d_0} & A\otimes A \ar[r] & 0,
}\]
where $V=\Bbbk\{x,y\}$ and $R=\Bbbk r$. The differentials are \\
$\begin{array}{ll}
 d_0(1\otimes v\otimes 1) \ = & v\otimes 1 - 1\otimes v,\ \forall v\in V,\\
 d_1(1\otimes r\otimes 1) \ = & y\otimes x\otimes 1 + 1\otimes y\otimes x - x\otimes y\otimes 1 - 1\otimes x\otimes y\\
 & \ +\ \frac{1}{2}\  x\otimes x\otimes 1 + \frac{1}{2}\ 1\otimes x\otimes x, 
\end{array}$\\ \\
and the quasi-isomorphism is the one induced by $\mu:A\otimes A\to A$.

We only need to define an action of $\Bbbk\Z=\Bbbk[t,t\inv]$ and a compatible internal $\Z$-grading on the $A$-bimodules $A\otimes R\otimes A$ and $A\otimes V\otimes A$. We will define the action of $t$ on the elements of the form $1\otimes v\otimes 1$ and extend it by $t\cdot (a\cdot m \cdot b) = (t\cdot a)\cdot (t\cdot m)\cdot (t\cdot b)$. Analogously, we will define the grading on these elements and extend it additively over the tensor product.
\begin{itemize}
    \item $t\cdot (1\otimes x\otimes 1) = 1\otimes x\otimes 1,\ t\cdot (1\otimes y\otimes 1) = 1\otimes x\otimes 1 + 1\otimes y\otimes 1$ and $1\otimes x\otimes 1,1\otimes y\otimes 1$ have internal degree $1$,
    \item $t\cdot (1\otimes r\otimes 1) = 1\otimes r\otimes 1$ and $1\otimes r\otimes 1$ has internal degree $2$.
\end{itemize}
It is straightworward to check that the differentials preserve the action and the internal grading.
\end{ex}

\begin{ex} \label{exresolsuperjordan}
For the super Jordan plane, we will use the following resolution of $A$ as $A$-bimodule, computed in \cite{recasolotar}:

\[\xymatrix{ \cdots \ar[r]^-{d_3} & A\otimes V_3\otimes A \ar[r]^-{d_2}& A\otimes V_2\otimes A \ar[r]^-{d_1} & A\otimes V_1\otimes A \ar[r]^-{d_0} & A\otimes A \ar[r] & 0,
}\]
where $V_1=\Bbbk\{x,y\}$ and $V_n=\Bbbk\{x^n, y^2 x^{n-1}\},\ \forall n\geq 2$. It has differentials \\
$\begin{array}{ll}
 d_0(1\otimes v\otimes 1)= v\otimes 1 - 1\otimes v,\ \forall v\in V,\\
 d_1(1\otimes x^2\otimes 1)= 1\otimes x\otimes x + x\otimes x\otimes 1,\\
 d_1 (1 \otimes y^2 x \otimes 1) = y^2\otimes x\otimes 1 + y\otimes y\otimes x + 1\otimes y\otimes yx - xy\otimes y\otimes 1 - x\otimes y\otimes y\\
 \qquad \qquad \qquad \qquad \quad -\ 1\otimes x\otimes y^2 - xy\otimes x\otimes 1 - x\otimes y\otimes x - 1\otimes x\otimes yx, \\
 d_n (1\otimes x^{n+1}\otimes 1) = x\otimes x^n\otimes 1 + (-1)^{n+1}\ 1\otimes x^n\otimes x, \\
d_n (1\otimes y^2 x^n\otimes 1) = y^2\otimes x^n\otimes 1 + (-1)^{n+1}\ 1\otimes y^2 x^{n-1}\otimes x - x\otimes y^2 x^{n-1}\otimes 1 \\
\qquad \qquad \qquad \qquad \quad -\ xy\otimes x^n\otimes 1 - 1\otimes x^n\otimes y^2 - 1\otimes x^n\otimes yx,\ \forall n\geq 2,
\end{array}$

and quasi-isomorphism induced by $\mu:A\otimes A\to A$ as in the previous example.

Now, we define a compatible $\Bbbk\Z$-action and internal $\Z$-grading on the elements $1\otimes v\otimes 1$ of the $A$-bimodules $A\otimes V_n\otimes A$:
\begin{itemize}
    \item $t\cdot (1\otimes x\otimes 1) = -1\otimes x\otimes 1,\ t\cdot (1\otimes y\otimes 1) = 1\otimes x\otimes 1 - 1\otimes y\otimes 1$ and $1\otimes x\otimes 1,1\otimes y\otimes 1$ have internal degree $1$,
    \item $t\cdot (1\otimes x^n\otimes 1) = (-1)^n 1\otimes x^n\otimes 1$, and $1\otimes x^n\otimes 1$ has internal degree $n$, $\forall n\geq 2$.
    \item $t\cdot (1\otimes y^2 x^{n-1}\otimes 1) = -1\otimes x^n \otimes y + (-1)^{n-1} (x-y)\otimes x^n \otimes 1 + 1\otimes y^2 x^{n-1} \otimes 1$, and $1\otimes y^2 x^{n-1} \otimes 1$ has internal degree $n+1$, $\forall n\leq 2$.
    \end{itemize}

As in the previous example, it is straightforward to check the needed compatibilities.
\end{ex}

\medskip

So far in this section, we have worked on a duoidal structure in the category $\AkGYDA$ analogous to the braided monoidal structure in $\kGYD$. We have also defined a coduoid up to homotopy in $\khCh\cc\AkGYDA\dd$ using a free resolution $P_\bullet\to A$, in a similar way as we had previously defined a braided comonoid up to homotopy in $\hCh\cc\kGYD\dd$ using the complex $S_\bullet$. 

Our next step will be to define a functor $\homAA(-,\Bbbk):\AkGYDA \to \kGYD$ in a similar way to the definition of $\hom(-,\Bbbk):\kGYD\to\kGYD$ in Example \ref{exhomintkGYD}. Under weaker finiteness conditions than in \cite{mpsw}, it will coincide at the $\Bbbk$-linear level with the functor $\HomAA(-,\Bbbk)$, sending a free resolution $P_\bullet\to A$ to the cochain complex used to compute Hochschild cohomology. We will extend this to the monoidal structures on them. The following definition is an analogue in the context of duoidal categories of the concept of braided lax monoidal functors in braided monoidal categories.

\begin{defi} \label{defidoublelax}
Let $(\CC,\diamond,I,\star,J,\zeta)$ and $(\CC',\diamond',I',\star',J',\zeta')$ be duoidal categories. A \emph{double lax monoidal functor} from $\CC$ to $\CC'$
is a 5-tuple $(F,\varphi,\varphi_0,\gamma,\gamma_0)$, where
\begin{itemize}
    \item $(F,\varphi,\varphi_0)$ is a lax monoidal functor from $(\CC,\diamond, I)$ to $(\CC',\diamond', I')$.
    \item $(F,\gamma,\gamma_0)$ is a lax monoidal functor from $(\CC,\star, J)$ to $(\CC',\star', J')$.
    \item The following diagrams commute:
\end{itemize}
\[
\xymatrix{
 (F(X)\star' F(Y))\diamond' (F(Z)\star' F(T)) \ar[rrr]^{\zeta'_{F(X),F(Y),F(Z),F(T)}} \ar[d]_{\gamma_{X,Y}\diamond'\gamma_{Z,T}} & & & (F(X)\diamond' F(Z))\star' (F(Y)\diamond' F(T)) \ar[d]^{\varphi_{X,Z}\star'\varphi_{Y,T}} \\
 F(X\star Y)\diamond' F(Z\star T) \ar[d]_{\varphi_{X\star Y, Z\star T}} & & & F(X\diamond Z)\star' F(Y\diamond T) \ar[d]^{\gamma_{X\diamond Z,Y\diamond T}} \\
 F((X\star Y)\diamond(Z\star T)) \ar[rrr]_{F(\zeta_{X,Y,Z,T})} & & & F((X\diamond Z)\star(Y\diamond T)),
}\]
\[
\xymatrix{
I' \ar[r]^-{\Delta_{I'}} \ar[d]_{\varphi_0} & I'\star'  I' \ar[dd]^{\varphi_0\star'\varphi_0} \\
F(I) \ar[d]_{F(\Delta_I)}   &  \\
F(I\star I) & F(I)\star' F(I) \ar[l]^{\gamma_{I,I}},
} \quad 
\xymatrix{
J'\diamond' J' \ar[r]^-{\mu_{J'}} \ar[dd]_{\gamma_0\diamond'\gamma_0}  & J \ar[d]^{\gamma_0} \\
 & F(J)      \\
F(J)\diamond' F(J) \ar[r]_-{\varphi_{J,J}}& F(J\diamond J), \ar[u]_{F(\mu_J)}
} \quad
\xymatrix{
I' \ar[r]^{\zeta'_0} \ar[d]_{\varphi_0} & J' \ar[d]^{\gamma_0}\\ F(I) \ar[r]_{F(\zeta_0)} & F(J). 
}
\]
\end{defi}

While lax monoidal functors preserve monoids and braided lax monoidal functors preserve commutative monoids, we recall a proposition which states that double lax monoidal functors preserve duoids.

\begin{prop} \label{propdbllaxsendsduoid} \cite[Corollary 6.58]{aguiar}
Given two duoidal categories $(\CC,\diamond,I,\star,J,\zeta)$ and $(\CC',\diamond',I',\star',J',\zeta')$, let $(F,\varphi,\varphi_0,\gamma,\gamma_0)$ be a double lax monoidal functor from $\CC$ to $\CC'$, and let $(A,\mu,\eta,\nu,\iota)$ be a duoid in $\CC$. The 5-tuple $(F(A),F(\mu)\circ\varphi_{A,A},F(\eta)\circ\varphi_0,F(\nu)\circ\gamma_{A,A},F(\iota)\circ\gamma_0)$ is a duoid in $\CC'$.
\end{prop} 

In order to define the functor mentioned before Definition \ref{defidoublelax}, let us recall the internal hom in the category $\kGYD$ from Examples \ref{exhomintVectG} and \ref{exhomintkGYD}. Componentwise, it is as follows:
\[
\hom(V,Y)_h = \prod_{s\in G} \Homk(V_s, Y_{hs}) = \Hom_{\Vect^G}(V,Y[h]),
\]
Where the $G$-graded vector space $Y[h]$ is defined componentwise by $Y[h]_s=Y_{hs}$.

The next definition resembles this internal hom.

\begin{defi}
Let $A$ be an augmented algebra and $M,N$ be two $A$-bimodules in the monoidal category $\kGYD$, where $G$ is an abelian group. The Yetter-Drinfeld module $\hom_{AA}(M,N)$ is defined componentwise by
\[
\hom_{AA}(M,N)_h =  \Hom_{_A\!(\Vect^G)\!_A}(M,N[h]) = \Hom_{AA}(M,N)\cap\hom(M,N)_h,
\]
with the following action of $\Bbbk G$:
\[
(g\cdot f)(x)=g\cdot f(g^{-1}\cdot x),\ \forall x\in M,\ \forall g\in G,\ \forall f\in\hom_{AA}(M,N)_h.
\]
\end{defi}

There are several details implicit in this definition:
\begin{itemize}

\item We are using the isomorphism $\kGYD \simeq \cc_{\Bbbk G}\Mod\dd^G$ of Example \ref{exHYD} to define a Yetter-Drinfeld module by a $G$-grading and a $\Bbbk G$-action on each component.

\item The aforementioned isomorphism is also used to define the object $N[h]$. It can be checked directly that $N[h]$ is an $A$-bimodule in $\kGYD$ if $N$ is so.

\item In Remark \ref{remaforgetstrong}, we noticed that the forgetful functors $\kGYD\longrightarrow\Vect^G\longrightarrow\Vect$ are strong monoidal, so $A$ can be regarded as a monoid and $M,N,N[h]$ as $A$-bimodules in any of these categories.

\item It is straightforward to check that if $f$ is in $\hom_{AA}(M,N)_h$, then $g\cdot f$ belongs to $\hom_{AA}(M,N)_h$ as well. Thus, $\hom_{AA}(M,N)_h$ is a $\Bbbk G$-submodule of $\hom(M,N)_h$ for every $h\in G$, and consequently $\homAA(M,N)$ is a $G$-graded $\Bbbk G$-submodule of $\hom(M,N)$.
\end{itemize}

It can also be easily checked that precomposing with $A$-bimodule morphisms in the first coordinate and postcomposing with $A$-bimodule morphisms in the second one preserves the $G$-grading and the $\Bbbk G$-action. This allows us to regard $\homAA(-,-)$ as a bifunctor from $\AkGYDA^\flop \times \AkGYDA$ to $\kGYD$.

\begin{rema}
Let $M$ be a Yetter-Drinfeld module over $\Bbbk G$ isomorphic to the free $A$-bimodule $A\otimes V\otimes A$ as a bimodule in $\Vect^G$. The adjunction between the functor $A\otimes - \otimes A$ and the forgetful functor from ${_A\!(\Vect^G)\!_A}$ to $\Vect^G$ gives a canonical isomorphism
\[
\Hom_{_A\!(\Vect^G)\!_A}(A\otimes V\otimes A,N[h]) \simeq \Hom_{\Vect^G}(V,N[h]).
\]
This implies that
\[
\homAA(A\otimes V\otimes A,N) \simeq \hom(V,N)
\]
as $G$-graded vector spaces.
If, in addition, $V$ is a finite dimensional vector space, then
\[
\homAA(A\otimes V\otimes A,N) \simeq \hom(V,N) \simeq \Homk(V,N) \simeq \HomAA(A\otimes V\otimes A,N).
\]
Since all these isomorphisms are natural, so is their composition.
\end{rema}

\begin{prop} \label{prophomAAHomAA}
Let $A$ be an augmented algebra in the monoidal category $\kGYD$ for an abelian group $G$. Suppose that there exists a chain complex $P_\bullet\to A$ in $\AkGYDA$ which is a resolution of $A$ in $\AModA$, and such that for every $n\in\N_0$ the $A$-bimodule $P_n$ is isomorphic as a $G$-graded vector space to $A\otimes V_n\otimes A$, with $V_n$ finite dimensional. There is an isomorphism of $\Z$-graded vector spaces
\[
\H^\bullet(A,\Bbbk) \simeq H(\homAA(P_\bullet,\Bbbk)).
\]
\end{prop}
\noindent\emph{Proof:} 
Since $V_n$ is finite dimensional, by the previous remark  $\homAA(A\otimes V_n \otimes A, \Bbbk)$ is isomorphic to $\HomAA(A\otimes V_n \otimes A, \Bbbk)$.

Using that $P_\bullet\to A$ is a projective resolution in $\Vect$ gives
\[
\H^\bullet(A,\Bbbk) = \Ext^\bullet_{AA}(A,\Bbbk) = H(\HomAA(P_\bullet,\Bbbk)) \simeq H(\homAA(P_\bullet,\Bbbk)).
\]
\qed

Now we know that, under some finiteness conditions, the functor $\homAA(-,\Bbbk)$ can be used to compute Hochschild cohomology. The following step is giving this functor a double lax monoidal structure relating the coproducts $\omega$ and $\delta$ of Theorem \ref{thmPcoduoid} to the cup product of $\H^\bullet(A,\Bbbk)$.

\begin{defi} \label{defihomduoidal}
Let $(A,\mu,\eta,\Delta,\varepsilon)$ be a bimonoid in $\kGYD$, and let $M,N\in \AkGYDA$. 
\begin{itemize}
\item The morphism
$\varphi_{M.N}: \homAA(M,\Bbbk) \otimes \homAA(N,\Bbbk) \longrightarrow \homAA(N\odot M, \Bbbk)$ is defined by $\varphi_{M,N}(f\otimes g)=g\odot f$.
\item The morphism $\varphi_0:\Bbbk \longrightarrow \homAA(\Bbbk,\Bbbk)$ is defined by $\varphi_0(1)=\id_\Bbbk$.
\item The morphism
$\gamma_{M.N}: \homAA(M,\Bbbk) \otimes \homAA(N,\Bbbk) \longrightarrow \homAA(N\otimes_A M, \Bbbk)$ is defined by $\gamma_{M,N}(f\otimes g)=g\otimes_A f$.
\item The morphism $\gamma_0:\Bbbk \longrightarrow \homAA(A,\Bbbk)$ is defined by $\gamma_0(1)=\varepsilon$.
\end{itemize} 
\end{defi}

\begin{thm}
Let $A$ be a bimonoid in $\kGYD$. The 5-tuple $\bigl(\homAA(-,\Bbbk), \varphi,\varphi_0, \gamma,\gamma_0 \bigr)$ is a double lax monoidal functor from $\cc \AkGYDA^\flop, \todot,\Bbbk, \totimesA,A, \cc\zeta^A\dd^{\todot,\totimesA} \dd$ to the duoidal category $\Bigl( \kGYD, \otimes, \Bbbk,, \otimes, \Bbbk, \zeta^c\Bigr)$.
\end{thm}
\noindent\emph{Proof:} First of all, notice that $\varphi$ and $\gamma$ are natural transformations due to the functoriality of the products $\odot, \otimes_A$ and $\otimes$. Associativity an unitality of both lax monoidal structures can be checked directly, as well as the axioms involving the respective comonoid and monoid structures on the first and second units. Let us pay attention to the relation between the functor $\homAA(-,\Bbbk)$ and the interchange laws $\cc\zeta^A\dd^{\todot,\totimesA}$ and $\zeta^c$. We need to prove the commutativity of the following diagram for every $M,N,K,L\in\AkGYDA$, where for $X\in\{M,N,K,L\}$ the object $\homAA(X,\Bbbk)$ is denoted by $X^\circ$.
\[
\xymatrix{
 M^\circ\otimes N^\circ \otimes K^\circ \otimes L^\circ \ar[rrr]^{\id_{_{M^\circ}}\otimes c_{_{N^\circ,K^\circ}}\otimes \id_{_{L^\circ}}} \ar[d]_{\gamma_{_{M,N}}\otimes\gamma_{_{K,L}}} & & & M^\circ\otimes K^\circ\otimes N^\circ\otimes L^\circ \ar[d]^{\varphi_{_{M,K}}\otimes\varphi_{_{N,L}}} \\
 (M\totimesA N)^\circ\otimes (K\totimesA L)^\circ \ar[d]_{\varphi_{_{M\totimesA N, K\totimesA L}}} & & & (M\todot K)^\circ\otimes (N\todot L)^\circ \ar[d]^{\gamma_{_{M\todot K,N\todot L}}} \\
 \bigl[(M\totimesA N)\todot(K\totimesA L)\bigr]^\circ \ar[rrr]_{\left(\left(\zeta^A\right)^{\todot,\totimesA}_{M,K,N,L}\right)^*} & & & \bigl[(M\todot K)\totimesA(N\todot L)\bigr]^\circ,
}\]
which translates into
\[
\xymatrix{
 M^\circ\otimes N^\circ\otimes K^\circ\otimes L^\circ \ar[rrr]^{\id_{_{M^\circ}}\otimes c_{_{N^\circ,K^\circ}}\otimes \id_{_{L^\circ}}} \ar[d]_{\gamma_{_{M,N}}\otimes\gamma_{_{K,L}}} & & & M^\circ\otimes K^\circ\otimes N^\circ\otimes L^\circ \ar[d]^{\varphi_{_{M,K}}\otimes\varphi_{_{N,L}}} \\
 (N\otimes_A M)^\circ\otimes (L\otimes_A K)^\circ \ar[d]_{\varphi_{_{N\otimes_A M, L\otimes_A K}}} & & & (K \odot M)^\circ\otimes (L \odot N)^\circ \ar[d]^{\gamma_{_{K \odot M,L \odot N}}} \\
 \bigl[(L\otimes_A K)\odot (N\otimes_A M)\bigr]^\circ \ar[rrr]_{\left(\zeta^A_{L,N,K,M}\right)^*} & & & \bigl[(L \odot N)\otimes_A (K \odot M)\bigr]^\circ,
}\]
Now, let $a,b,c,d$ be elements of $G$. We will chase an element $f\otimes f'\otimes g\otimes g'$ along the diagram above, where $f\in\homAA(M_a,\Bbbk)$, $f'\in\homAA(N_b,\Bbbk)$, $g\in\homAA(K_c,\Bbbk)$ and $g'\in\homAA(L_d,\Bbbk)$ respectively. Notice that, since $\Bbbk$ is concentrated in degree $1_G$, the respective internal degrees of $f, f', g, g'$ are $a\inv, b\inv, c\inv, d\inv$.
Following the top and right sides of the diagram gives the element $(g'\odot f') \otimes_A ( b\inv\cdot g\odot f)$, while following its left and bottom sides gives $\cc (g'\otimes_A g)\odot (f'\otimes_A f) \dd \circ \zeta^A_{L,N,K,M}$. To show that both terms are equal, it suffices to evaluate them in a generic element $(w'\odot v')\otimes_A (w\odot v)$, where $v,v',w,w'$ lie in $M_a,N_b,K_c,L_d$ respectively. Evaluating the first term yields
\[
\begin{array}{rcl}
(g'\odot f') \otimes_A ( b\inv\cdot g\odot f) \cc (w'\odot v')\otimes_A (w\odot v) \dd  &=& g'(w')f'(v')(b\inv\cdot g)(w)f(v) \\
&=& g'(w')f'(v')g(b\cdot w)f(v),
\end{array}
\]
since the action of $\Bbbk G$ on $\Bbbk$ is trivial, while evaluating the second term yields
\begin{multline*}
\cc (g'\otimes_A g)\odot (f'\otimes_A f)\dd  \cc \zeta^A_{L,N,K,M} ((w'\otimes_A v')\odot (w\otimes v)) \dd \\
=\ \cc (g'\otimes_A g)\odot (f'\otimes_A f)\dd \cc (w'\otimes_A b\cdot w)\odot (v'\otimes_A v) \dd\ 
=\ g'(w')g(b\cdot w)f'(v')f(v).
\end{multline*}
This shows that both compositions of the diagram coincide when evaluated in a generic element $f\otimes f'\otimes g\otimes g'$, proving the compatibility of the functor $\homAA(-,\Bbbk)$ with the interchange laws.
\qed

\medskip

In the following theorem, we prove the result on graded braided commutativity up to $\Bbbk$-linear homotopy which leads to the desired property of graded braided commutativity of the Hochschild cohomology with trivial coefficients, under hypotheses that include Nichols algebras such as the Jordan and super Jordan plane. 

\begin{thm} \label{thmmain}
Let $A$ be a bimonoid in the braided monoidal category $\left(\kGYD, \otimes, \Bbbk, c\right)$ for an abelian group $G$. Suppose that there exists a complex $P_\bullet\to A$ in $\AkGYDA$ which is a projective resolution of $A$ in $\AModA$, and such that for every $n\in\N_0$ the $A$-bimodule $P_n$ is isomorphic as a $G$-graded vector space to $A\otimes V_n\otimes A$, where $V_n$ is finite dimensional. Then the differential graded algebra $\HomAA(P_\bullet,\Bbbk)$ with the opposite of the cup product is graded braided commutative up to $\Bbbk$-linear homotopy.
\end{thm}
\noindent\emph{Proof:}
In the same way as in Theorem \ref{thmMPSWhomotopy}, one can consider the opposite of the category $\khCh\cc\AkGYDA\dd$, which is isomorphic to the category of cochain complexes in $\AkGYDA$ with cochain morphisms in $\AModA$ in the opposite direction. 

On the other hand, let $\khCh\cc\kGYD\dd$ denote the category whose objects are cochain complexes of Yetter-Drinfeld modules over $\Bbbk G$, and whose morphisms are $\Bbbk$-linear homotopy classes of $\Bbbk$-linear cochain morphisms. It is a braided monoidal category with the tensor product defined as in Remark \ref{remaChCochCC} and the graded braiding.

It is posible to define componentwise a cochain version of the functor $\homAA(-,\Bbbk)$ from the duoidal category $\cc \khCh\cc\AkGYDA\dd^\flop, \todot,\Bbbk_0, \totimesA,A_0, \cc\zeta^A\dd^{\todot,\totimesA, gr} \dd$ to the duoidal category $\Bigl( \khCoch\cc\kGYD\dd, \otimes, \Bbbk_0,\otimes, \Bbbk_0, \cc\zeta^c\dd^{gr}\Bigr)$. All the axioms of a double lax monoidal functor hold at the cochain level because they do componentwise, and the graded interchange laws are compatible with the signs of the differentials given as in Remark \ref{remaChCochCC}, observing that they are being taken on the transpose product.

Translating the result of Theorem \ref{thmPcoduoid} to this context, one has that $\cc P_\bullet, \delta, f_0\circ\varepsilon_0, \omega, f_0 \dd$ is a duoid in $\khCoch\cc \AkGYDA^\flop, \todot,\Bbbk_0, \totimesA,A_0, \cc\zeta^A\dd^{\todot,\totimesA,gr} \dd$. Proposition \ref{propdbllaxsendsduoid} implies that the double lax functor $\homAA(-,\Bbbk)$ sends the duoid $P_\bullet$ to a duoid $\homAA(P_\bullet,\Bbbk)$ in $\Bigl(\underline{\Coch}\cc\kGYD\dd, \otimes, \Bbbk_0,\otimes, \Bbbk_0, \cc\zeta^c\dd^{gr}\Bigr)$. This duoid turns out to be a commutative monoid in the braided monoidal category $\Bigl(\underline{\Coch}\cc\kGYD\dd, \otimes, \Bbbk_0, c^{gr}\Bigr)$ by the Eckmann-Hilton argument in Example \ref{exEckmannHilton}. As proved in Proposition \ref{prophomAAHomAA}, the finiteness conditions on the free resolution $P_\bullet$ imply that $\homAA(P_\bullet,\Bbbk)$ is isomorphic to $\HomAA(P_\bullet,\Bbbk)$ as a cochain complex of $\Bbbk$-vector spaces. Regarding the monoid structure on $\HomAA(P_\bullet,\Bbbk)$ as the one induced by the comonoid $(P_\bullet,\omega,(\id_A)_0)$ via the lax monoidal functor $(\homAA(-,\Bbbk),\gamma,\gamma_0)$, one deduces that it coincides with the opposite of the usual cup product as defined in Example \ref{exomega}.
\qed

\begin{cor} \label{cormain}
Under the conditions above, $\H^\bullet(A,\Bbbk)$ is a graded braided commutative algebra in $\kGYD$.
\end{cor}

The projective resolutions of the Jordan and super Jordan plane shown respectively in Examples \ref{exresoljordan} and \ref{exresolsuperjordan} verify the hypotheses of the previous theorem, therefore their respective Hochschild cohomology algebras with the opposite of the cup product are graded braided commutative algebras in $\kZYD$.

\appendix

\section{Computations for the Jordan and super Jordan planes}

In this Appendix we make the computations needed to check the commutativity of the second diagram in Remark \ref{remacheckcomm}, that is of:

\[
\xymatrix{
\HomAA(A^{\otimes p+2},\Bbbk)\otimes\HomAA(A^{\otimes q+2},\Bbbk)  \ar[r]^-{\smallsmile_{p,q}} \ar[dr]_{(-1)^{pq}\smallsmile_{p,q}} &
\HomAA(A^{\otimes p+q+2},\Bbbk) \ar[d]^{\HomAA\cc \id_{_A}\otimes c_{_{A^{\otimes p},A^{\otimes q}}}\otimes\id_{_A},\Bbbk\dd} \\ & \HomAA(A^{\otimes p+q+2},\Bbbk)..}
\]
In order to check the commutativity of this diagram after taking cohomology, we will choose in each case a suitable projective resolution of $A$ as $A$-bimodule, denoted $P_\bullet\to A$, along with comparison morphisms $f_\bullet:P_\bullet\to B_\bullet(A)$ and $g_\bullet:B_\bullet(A)\to P_\bullet$.

For any two cocycles $\alpha$ and $\beta$, let us write $\alpha\sim\beta$ if their cohomology classes coincide. Commutativity means that
\[
 \HomAA\cc \id_{_A}\otimes c_{_{A^{\otimes p},A^{\otimes q}}}\otimes\id_{_A},\Bbbk\dd \cc \smallsmile_{p,q}(\psi'\otimes\varphi') \dd \sim  (-1)^{pq}\smallsmile_{p,q}(\psi'\otimes\varphi'),
\]
 for every $\psi'\in\HomAA(A^{\otimes p+2},\Bbbk)$ and $\varphi'\in\HomAA(A^{\otimes q+2},\Bbbk)$, that is:

\[
  (\psi'\smallsmile\varphi') \circ \cc \id_{_A}\otimes c_{_{A^{\otimes p},A^{\otimes q}}}\otimes\id_{_A} \dd \sim (-1)^{pq} (\psi'\smallsmile\varphi').
\]
Replacing the bar resolution by $P_\bullet$ this reads as
\[
 (\psi\circ g_p\smallsmile\varphi\circ g_q) \circ \cc \id_{_A}\otimes c_{_{A^{\otimes p},A^{\otimes q}}}\otimes\id_{_A} \dd \circ f_{p+q} \sim 
(-1)^{pq} (\psi\circ g_p\smallsmile\varphi\circ g_q) \circ f_{p+q},
\]
for every $\psi\in\HomAA(P_p,\Bbbk)$ and $\varphi\in\HomAA(P_q,\Bbbk)$.
This equivalence of $A$-bimodule morphisms from $P_{p+q}$ to $\Bbbk$ is what we will check.
\bigskip

Let us start with the Jordan plane. We will use the resolution of 
Example \ref{exresoljordan}:
\[\xymatrix{
0 \ar[r] & A\otimes R\otimes A \ar[r]^-{d_1} & A\otimes V\otimes A \ar[r]^-{d_0} & A\otimes A \ar[r] & 0,
}\]
where $V=\Bbbk\{x,y\}$ and $R=\Bbbk r$.

It is possible to define comparison morphisms $f_\bullet:P_\bullet \to B_\bullet(A)$ and $g_\bullet:B_\bullet(A) \to P_\bullet$ such that:\\

$\begin{array}{l}
 f_0 = g_0 = \id_{A\otimes A},\\
 f_1(1\otimes v\otimes 1) = 1\otimes v\otimes 1,\ \forall v\in V,\\
 g_1(1\otimes x^k y^l\otimes 1) = \sum_{i=0}^{k-1} x^i\otimes x\otimes x^{k-1-i} y^l + \sum_{j=0}^{l-1} x^k \otimes y\otimes y^{l-1-j}\\
 f_2(1\otimes r\otimes 1) = 1\otimes y\otimes x\otimes 1 - 1\otimes x\otimes y\otimes 1 + \frac{1}{2}\  1\otimes x\otimes x\otimes 1.
\end{array}$

\medskip

Notice that the commutativity needs to be checked only for the cup products of $\H^1(A,\Bbbk)$ with itself, since $\H^0(A,\Bbbk)$ consists only of scalar multiples of the unit object and the product of an element in $\H^2(A,\Bbbk)$ with one in $\H^1(A,\Bbbk)$ or $\H^2(A,\Bbbk)$ is zero since $\H^3(A,\Bbbk)=0$. Now, let $\psi,\varphi\in\HomAA(P_1,\Bbbk)$ be two 1-cocycles:
\[
\begin{array}{l}
(\psi\circ g_1\smallsmile\varphi\circ g_1) (f_2(1\otimes r\otimes 1))\\
\qquad =\  (\psi\circ g_1\smallsmile\varphi\circ g_1) \bigl( 1\otimes y\otimes x\otimes 1  - 1\otimes x\otimes y\otimes 1 +  \frac{1}{2}\  1\otimes x\otimes x\otimes 1 \bigr) \\
\qquad =\ \varphi(g_1(1\otimes y\otimes 1))\psi(g_1(1\otimes x\otimes 1)) - \varphi(g_1(1\otimes x\otimes 1))\psi(g_1(1\otimes y\otimes 1)) \\
\qquad \qquad +\ \frac{1}{2}\  \varphi(g_1(1\otimes x\otimes 1))\psi(g_1(1\otimes x\otimes 1)) \\
\qquad =\ \varphi(1\otimes y\otimes 1)\psi(1\otimes x\otimes 1) - \varphi(1\otimes x\otimes 1)\psi(1\otimes y\otimes 1) \\
\qquad \qquad +\ \frac{1}{2}\  \varphi(1\otimes x\otimes 1)\psi(1\otimes x\otimes 1).
\end{array}
\]
To compute the element we want to compare with the one above, recall that we write $\Bbbk\Z=\Bbbk[t,t\inv]$, and the braiding is given by $c_{_{V,W}}(v\otimes w)=t^n\cdot w\otimes v$, where $n$ is the degree of $v$. From now on, let $c_{p,q}$ denote the map $\id_{_A} \otimes c_{_{A^{\otimes p}, A^{\otimes q}}} \otimes \id_{_A}$ for every $p,q\in\N_0$.
\[
\begin{array}{l}
 (\psi\circ g_1\smallsmile\varphi\circ g_1) (c_{1,1} (f_2(1\otimes r\otimes 1)) \\
 \qquad =\  (\psi\circ g_1\smallsmile\varphi\circ g_1) \bigl( c_{1,1} \bigl( 1\otimes y\otimes x\otimes 1 - 1\otimes x\otimes y\otimes 1 \ +\  \frac{1}{2}\  1\otimes x\otimes x\otimes 1 \bigr)\bigr) \\
\qquad =\ (\psi\circ g_1\smallsmile\varphi\circ g_1) \bigl( 1\otimes t\cdot x\otimes y\otimes 1 - 1\otimes t\cdot y \otimes x \otimes 1 \ +\ \frac{1}{2}\ 1\otimes t\cdot x\otimes x\otimes 1 \bigr) \\
\qquad =\ (\varphi\circ g_1\smallsmile\psi\circ g_1) \bigl( 1\otimes x\otimes y\otimes 1 - 1\otimes (x+y) \otimes x \otimes 1 \ +\ \frac{1}{2}\ 1\otimes x\otimes x\otimes 1 \bigr) \\
\qquad  =\ -\varphi(1\otimes y\otimes 1)\psi(1\otimes x\otimes 1) + \varphi(1\otimes x\otimes 1)\psi(1\otimes y\otimes 1) \\
\qquad \qquad -\ \frac{1}{2}\  \varphi(1\otimes x\otimes 1)\psi(1\otimes x\otimes 1)\\
\qquad =\ -\ (\psi\circ g_1\smallsmile\varphi\circ g_1) (f_2(1\otimes r\otimes 1))
\end{array}
\]
as we wanted to prove.

\bigskip

The computation for the super Jordan plane is more involved, since the minimal resolution, which is free, is as follows -see Example \ref{exresolsuperjordan}-:
\[\xymatrix{ \cdots \ar[r]^-{d_3} & A\otimes V_3\otimes A \ar[r]^-{d_2}& A\otimes V_2\otimes A \ar[r]^-{d_1} & A\otimes V_1\otimes A \ar[r]^-{d_0} & A\otimes A \ar[r] & 0,
}\]
where $V_1=\Bbbk\{x,y\}$ and $V_n=\Bbbk\{x^n, y^2 x^{n-1}\},\ \forall n\geq 2$. 

For checking the commutativity condition, we will use the comparison morphisms $f_\bullet:P_\bullet\to B_\bullet(A)$ and $g_\bullet:B_\bullet(A)\to P_\bullet$ defined in \cite[Section 5]{recasolotar}. Since the chain morphism $g_\bullet$ is only partially defined in \cite{recasolotar}, we need to extend the definition of $g_n$ for $n\geq 2$ as follows, :
\[\begin{array}{l}
g_n(1\otimes x^{\otimes i}\otimes xy\otimes x^{\otimes n-1-i}\otimes 1) = 0 \quad \forall i:\ 0\leq i\leq n-2,\\
g_n(1\otimes x^{\otimes n-1}\otimes xy\otimes 1) = 1\otimes x^n\otimes y.
\end{array}
\]
Again, $\H^0(A,\Bbbk)= \Bbbk$, so we will check the condition for cocycles in degrees greater or equal to $1$. The cup products to be computed lie in $\H^n(A,\Bbbk)$ for $n\geq 2$, so we are going to evaluate them on the basis of $P_n$. Let $\psi\in\HomAA(P_p,\Bbbk)$ and $\varphi\in\HomAA(P_q,\Bbbk)$ be two cocycles, for $p,q\geq 1$. 

For $p=q=1$:

Evaluating in $1\otimes x^2\otimes 1$ gives
\[
\begin{array}{rcl}
(\psi\circ g_1 \smallsmile \varphi\circ g_1)(f_2(1\otimes x^2 \otimes 1)) &=& (\psi\circ g_1 \smallsmile \varphi\circ g_1)(1\otimes x^{\otimes 2}\otimes 1) \\
&=& \varphi(g_1(1\otimes x\otimes 1)) \psi(g_1(1\otimes x\otimes 1)) \\
&=& \varphi(1\otimes x\otimes 1) \psi(1\otimes x\otimes 1),
\end{array}
\]
\[
\begin{array}{rcl}
(\psi\circ g_1 \smallsmile \varphi\circ g_1)(c_{1,1}(f_2(1\otimes x^2\otimes 1))) 
&=& (\psi\circ g_1 \smallsmile \varphi\circ g_1)(c_{1,1}(1\otimes x\otimes x\otimes 1))\\
&=& (\psi\circ g_1 \smallsmile \varphi\circ g_1)(1\otimes t\cdot x\otimes x\otimes 1)\\
&=& (\psi\circ g_1 \smallsmile \varphi\circ g_1)(1\otimes (-x)\otimes x\otimes 1)\\
&=& -\ \varphi(g_1(1\otimes x\otimes 1)) \psi(g_1(1\otimes x\otimes 1)) \\
&=& -\ \varphi(1\otimes x\otimes 1) \psi(1\otimes x\otimes 1) \\
&=& -\ (\psi\circ g_1 \smallsmile \varphi\circ g_1)(f_2(1\otimes x^2 \otimes 1)),
\end{array}
\]
while evaluating in $1\otimes y^2 x\otimes 1$ gives
\[\begin{array}{l}
(\psi\circ g_1 \smallsmile \varphi\circ g_1)(f_2(1\otimes y^2 x\otimes 1)) \\
\qquad =\ (\psi\circ g_1 \smallsmile \varphi\circ g_1) (1\otimes y\otimes yx\otimes 1 - 1\otimes x\otimes y^2 \otimes 1 - 1\otimes x\otimes yx\otimes 1) \\
\qquad =\ \varphi(g_1(1\otimes y\otimes 1)) \psi(g_1(1\otimes yx\otimes 1)) - \varphi(g_1(1\otimes x\otimes 1)) \psi(g_1(1\otimes y^2\otimes 1)) \\
\qquad \qquad -\ \varphi(g_1(1\otimes x\otimes 1)) \psi(g_1(1\otimes yx\otimes 1))\\
\qquad =\ 0,
\end{array}
\]
since $g_1(1\otimes a\otimes 1)=0$ if $a$ is homogeneous of degree greater than $1$. On the other hand,
\[\begin{array}{l}
(\psi\circ g_1 \smallsmile \varphi\circ g_1)(c_{1,1}(f_2(1\otimes y^2 x\otimes 1))) \\
  =\ (\psi\circ g_1 \smallsmile \varphi\circ g_1)(c_{1,1}(1\otimes y\otimes yx\otimes 1 - 1\otimes x\otimes y^2 \otimes 1 - 1\otimes x\otimes yx\otimes 1))\\
  =\ (\psi\circ g_1 \smallsmile \varphi\circ g_1) (1\otimes t\cdot yx \otimes y\otimes 1 - 1\otimes t\cdot y^2 \otimes x\otimes 1 - 1\otimes t\cdot yx\otimes x\otimes 1)\\
  =\ (\psi\circ g_1 \smallsmile \varphi\circ g_1) (1\otimes yx \otimes y\otimes 1 - 1\otimes (y^2-xy-yx) \otimes x\otimes 1 - 1\otimes yx\otimes x\otimes 1) \\
  = \varphi(g_1(1\otimes yx\otimes 1)) \psi(g_1(1\otimes y\otimes 1) - \varphi(g_1(1\otimes (y^2-xy-yx) \otimes 1)) \psi(g_1(1\otimes x\otimes 1) \\
  \qquad -\ \varphi(g_1(1\otimes yx\otimes 1)) \psi(g_1(1\otimes x\otimes 1))\\
  =\ 0\ =\ -\ (\psi\circ g_1 \smallsmile \varphi\circ g_1)(f_2(1\otimes y^2 x\otimes 1)).
\end{array}
\]
For $p=1,q\geq 2$, evaluating in $1\otimes x^{1+q}\otimes 1$ gives on one hand
\[
\begin{array}{rcl}
  (\psi \circ g_1 \smallsmile \varphi \circ g_q) (f_{1+q} (1\otimes x^{1+q}\otimes 1)) &=& (\psi \circ g_1 \smallsmile \varphi \circ g_q) (1\otimes x^{\otimes 1+q}\otimes 1)\\
  &=& \varphi (g_q(1\otimes x^{\otimes q}\otimes 1)( \psi (g_1(1\otimes x\otimes 1))\\
  &=& \varphi(1\otimes x^{\otimes q}\otimes 1) \psi(1\otimes  1),
\end{array}
\]
and on the other hand
\[
\begin{array}{l}
(\psi \circ g_1 \smallsmile \varphi \circ g_q) (c_{1,q} (f_{1+q} (1\otimes x^{1+q}\otimes 1))) \\
 \qquad \qquad \qquad \qquad \qquad \qquad \qquad \qquad =\ (\varphi \circ g_1 \smallsmile \psi \circ g_q) (c_{1,q} (1\otimes x^{\otimes 1+q}\otimes 1))\\
 \qquad \qquad \qquad \qquad \qquad \qquad \qquad \qquad =\ (\varphi \circ g_q \smallsmile \psi \circ g_1) (c_{1,q} (1\otimes x\otimes x^{\otimes q}\otimes 1))\\
 \qquad \qquad \qquad \qquad \qquad \qquad \qquad \qquad =\  (\psi \circ g_1 \smallsmile \varphi \circ g_q) (1\otimes t\cdot (x^{\otimes q})\otimes x\otimes 1)\\
 \qquad \qquad \qquad \qquad \qquad \qquad \qquad \qquad =\ (\psi \circ g_1 \smallsmile \varphi \circ g_q) (1\otimes (t\cdot x)^{\otimes q}\otimes x\otimes 1)\\
 \qquad \qquad \qquad \qquad \qquad \qquad \qquad \qquad =\ (\psi \circ g_1 \smallsmile \varphi \circ g_q) (1\otimes (-x)^{\otimes q}\otimes x\otimes 1) \\
 \qquad \qquad \qquad \qquad \qquad \qquad \qquad \qquad =\ (-1)^q\ (\psi \circ g_1 \smallsmile \varphi \circ g_q) (1\otimes x^{\otimes 1+q}\otimes 1)\\
 \qquad \qquad \qquad \qquad \qquad \qquad \qquad \qquad =\ (-1)^q\ (\psi \circ g_1 \smallsmile \varphi \circ g_q) (f_{1+q} (1\otimes x^{1+q}\otimes 1)),
\end{array}
\]
while evaluating in $1\otimes y^2 x^q\otimes 1$ gives
\[\begin{array}{l}
(\psi g_1\smallsmile \varphi g_q) (f_{1+q} (1\otimes y^2 x^q\otimes 1)) \\ \qquad
  =\ (\psi g_1\smallsmile \varphi g_q)(1\otimes y\otimes yx\otimes x^{\otimes q-1}\otimes 1 - 1\otimes x\otimes (y^2+yx)\otimes x^{\otimes q-1}\otimes 1  \\ \qquad
   \qquad +\ \sum_{i=0}^{q-2}(-1)^i\ 1\otimes x^{\otimes 2+i}\otimes(y^2+yx)\otimes x^{\otimes q-2-i}\otimes 1) \\ \qquad
  =\ \varphi(g_q(1\otimes y\otimes yx\otimes x^{\otimes q-2}\otimes 1)) \psi(g_1(1\otimes x \otimes 1))  \\ \qquad
  \qquad -\ \varphi(g_q(1\otimes x\otimes (y^2+yx)\otimes x^{\otimes q-2}\otimes 1)) \psi(g_1(1\otimes x \otimes 1))  \\ \qquad
  \qquad +\ \sum_{i=0}^{q-3} (-1)^i\ \varphi(g_q(1\otimes x^{\otimes 2+i}\otimes(y^2+yx)\otimes x^{\otimes q-3-i}\otimes 1)) \psi(g_1(1\otimes x\otimes 1))  \\ \qquad
  \qquad +\ (-1)^{q-2}\ \varphi(g_q(1\otimes x^{\otimes q}\otimes 1)) \psi(g_1(1 \otimes(y^2+yx)\otimes 1)) \\ \qquad
  =\ \varphi(1\otimes y^2 x^{q-1}\otimes 1)\psi(1\otimes x\otimes 1),
\end{array}
\]
\[
\begin{array}{l}
(\psi g_1\smallsmile \varphi g_q) (c_{1,q} (f_{1+q} (1\otimes y^2x^q\otimes 1))) \\ 
\qquad =\ (\psi g_1\smallsmile \varphi g_q) (c_{1,q}(1\otimes y\otimes yx\otimes x^{\otimes q-1}\otimes 1- 1\otimes x\otimes (y^2+yx)\otimes x^{\otimes q-1}\otimes 1 \\
 \qquad \qquad +\ \sum_{i=0}^{q-2}(-1)^i\ 1\otimes x^{\otimes 2+i}\otimes(y^2+yx)\otimes x^{\otimes q-2-i}\otimes 1))\\
\qquad =\ (\psi g_1\smallsmile \varphi g_q) (1\otimes t\cdot yx\otimes (t\cdot x)^{\otimes q-1}\otimes y\otimes 1 \\
 \qquad \qquad -\ 1\otimes  t\cdot (y^2+yx)\otimes (t\cdot x)^{\otimes q-1}\otimes x \otimes 1 \\
 \qquad \qquad +\ \sum_{i=0}^{q-2}(-1)^i\ 1\otimes (t\cdot x)^{\otimes 1+i}\otimes t\cdot (y^2+yx)\otimes (t\cdot x)^{\otimes q-2-i}\otimes x\otimes 1)\\
\qquad  =\ (\psi g_1\smallsmile \varphi g_q) (1\otimes yx\otimes (-x)^{\otimes q-1}\otimes y\otimes 1 \\
 \qquad \qquad -\ 1\otimes (y^2-xy)\otimes (-x)^{\otimes q-1}\otimes x \otimes 1 \\
 \qquad \qquad +\ \sum_{i=0}^{q-2}(-1)^i\ 1\otimes (-x)^{\otimes 1+i}\otimes(y^2-xy)\otimes (-x)^{\otimes q-2-i}\otimes x\otimes 1)\\
\qquad  =\ (-1)^{q-1}\ \varphi(g_q(1\otimes yx\otimes x^{\otimes q-1})) \psi(g_1(1\otimes y\otimes 1)) \\
 \qquad \qquad -\ (-1)^{q-1}\ \varphi(g_q(1\otimes y^2\otimes x^{\otimes q-1}\otimes 1)) \psi(g_1(1\otimes 1)) \\
 \qquad \qquad +\ \varphi(g_q(1\otimes x^{\otimes q-1}\otimes xy\otimes 1)) \psi(g_1(1\otimes x\otimes 1)) \\
\qquad  =\ (-1)^q\  \varphi(1\otimes y^2 x^{q-1}\otimes 1) \psi(1\otimes x\otimes 1)\\
\qquad =\ (-1)^q\ (\psi g_1\smallsmile \varphi g_q) (f_{1+q} (1\otimes y^2x^q\otimes 1)) 
\end{array}
\]
For $p\geq 2,q=1$, evaluating in $1\otimes x^{p+1}\otimes 1$ gives
\[
\begin{array}{rcl}
  (\psi \circ g_p \smallsmile \varphi \circ g_1) (f_{p+1} (1\otimes x^{p+1}\otimes 1)) &=& (\psi \circ g_p \smallsmile \varphi \circ g_1) (1\otimes x^{\otimes p+1}\otimes 1)\\
  &=&  \varphi (g_1(1\otimes x^{\otimes 1}\otimes 1)( \psi (g_p(1\otimes x^{\otimes p}\otimes 1))\\
  &=& \varphi(1\otimes x\otimes 1) \psi(1\otimes x^p\otimes 1)
  \end{array}
  \]
  \[
  \begin{array}{l}
  (\psi \circ g_p \smallsmile \varphi \circ g_1) (c_{p,1} (f_{p+1} (1\otimes x^{p+1}\otimes 1)))\ =\ (\psi \circ g_p \smallsmile \varphi \circ g_1) (c_{p,1} (1\otimes x^{\otimes p+1}\otimes 1))\\
  \qquad \qquad \qquad \qquad \qquad \qquad \qquad \qquad =\ (\psi \circ g_p \smallsmile \varphi \circ g_1) (c_{p,1} (1\otimes x^{\otimes p}\otimes x\otimes 1))\\
  \qquad \qquad \qquad \qquad \qquad \qquad \qquad \qquad =\  (\psi \circ g_p \smallsmile \varphi \circ g_1) (1\otimes t^p\cdot x\otimes x^{\otimes p}\otimes 1)\\
  \qquad \qquad \qquad \qquad \qquad \qquad \qquad \qquad =\ (\psi \circ g_p \smallsmile \varphi \circ g_1) (1\otimes (-1)^p x\otimes x^{\otimes p}\otimes 1) \\
  \qquad \qquad \qquad \qquad \qquad \qquad \qquad \qquad =\ (-1)^{p}(\psi \circ g_p \smallsmile \varphi \circ g_1) (1\otimes x^{\otimes p+1}\otimes 1)\\
  \qquad \qquad \qquad \qquad \qquad \qquad \qquad \qquad =\ (-1)^{1p} (\psi \circ g_p \smallsmile \varphi \circ g_1) (f_{p+1} (1\otimes x^{p+1}\otimes 1)).
 \end{array}
\]
Evaluating in $1\otimes y^2 x^p\otimes 1$ gives
\[\begin{array}{l}
(\psi g_p\smallsmile \varphi g_1) (f_{p+1} (1\otimes y^2 x^p\otimes 1))\\ \qquad \qquad
=\ (\psi g_p\smallsmile \varphi g_q)(1\otimes y\otimes yx\otimes x^{\otimes p-1}\otimes 1 \\ \qquad \qquad
 \qquad \qquad -\ 1\otimes x\otimes (y^2+yx)\otimes x^{\otimes p-1}\otimes 1 \\ \qquad \qquad
 \qquad \qquad +\ \sum_{i=0}^{p-2}(-1)^i\ 1\otimes x^{\otimes 2+i}\otimes(y^2+yx)\otimes x^{\otimes p-2-i}\otimes 1)\\ \qquad \qquad
-\ \varphi(g_1(1\otimes y\otimes 1)) \psi(g_p(1\otimes yx\otimes x^{\otimes p-1}\otimes 1))\\ \qquad \qquad
 \qquad -\  \varphi(g_1(1\otimes x\otimes\ 1)) \psi(g_p(1 \otimes (y^2+yx)\otimes x^{\otimes p-1}\otimes 1))\\ \qquad \qquad
 \qquad +\ \sum_{i=0}^{p-2} (-1)^i\ \varphi(g_1(1\otimes x\otimes 1)) \psi(g_p(1\otimes x^{\otimes 1+i}\otimes(y^2+yx)\otimes x^{\otimes p-2-i}\otimes 1)).\\ \qquad \qquad
 =\ - \varphi(1\otimes y^2 x^{p-1}\otimes 1) \psi(1\otimes y^2 x^{p-1}\otimes 1),
  (\psi g_p\smallsmile \varphi g_1) (c_{p,1}(f_{p+1}(1\otimes y^2 x^p\otimes 1))) \\ \qquad
 =\ (\psi g_p\smallsmile \varphi g_1) (c_{p,1}(1\otimes y\otimes yx\otimes x^{\otimes p-1}\otimes 1 \\ \qquad
 \qquad \qquad -\ 1\otimes x\otimes (y^2+yx)\otimes x^{\otimes p-1}\otimes 1 \\ \qquad
 \qquad \qquad +\ \sum_{i=0}^{p-2}(-1)^i\ 1\otimes x^{\otimes 2+i}\otimes(y^2+yx)\otimes x^{\otimes p-2-i}\otimes 1))\\ \qquad
 =\ (\psi g_p\smallsmile \varphi g_1)(1\otimes t^{p+1}\cdot x\otimes  y\otimes yx\otimes x^{\otimes p-2}\otimes 1 \\ \qquad
 \qquad \qquad -\ 1\otimes t^{p+1}\cdot x\otimes x\otimes (y^2+yx)\otimes x^{\otimes p-2}\otimes 1 \\ \qquad
 \qquad \qquad +\ \sum_{i=0}^{p-3}(-1)^i\ 1\otimes t^{p+1}\cdot x\otimes x^{\otimes 2+i}\otimes(y^2+yx)\otimes x^{\otimes p-3-i}\otimes 1\\ \qquad
 \qquad +\ (-1)^{p-2}\ 1\otimes t^p\cdot (y^2 + xy) \otimes x^{\otimes p}\otimes 1) 
 \end{array}
\]
\[
\begin{array}{l}
 \qquad=\ (\psi g_p\smallsmile \varphi g_1)(1\otimes (-1)^{p+1}\ x\otimes  y\otimes yx\otimes x^{\otimes p-2}\otimes 1 \\ \qquad
 \qquad \qquad -\ 1\otimes (-1)^{p+1}\ x\otimes x\otimes (y^2+yx)\otimes x^{\otimes p-2}\otimes 1 \\ \qquad
 \qquad \qquad +\ \sum_{i=0}^{p-3}(-1)^i\ 1\otimes (-1)^{p+1}\ x\otimes x^{\otimes 2+i}\otimes(y^2+yx)\otimes x^{\otimes p-3-i}\otimes 1\\ \qquad
 \qquad +\ (-1)^{p-2}\ 1\otimes (y^2 -(p-1)yx-pxy) \otimes x^{\otimes p}\otimes 1) \\ 
 \end{array}
\]
\[
\begin{array}{l}
 \qquad
 =\ (-1)^{p+1}   \varphi(g_1(1\otimes x\otimes 1)) \psi(g_p(1\otimes y\otimes yx\otimes x^{\otimes p-2}\otimes 1)) \\ \qquad
 =\ (-1)^{p+1} \varphi(1\otimes y^2 x^{p-1}\otimes 1) \psi(1\otimes y^2 x^{p-1}\otimes 1)\\
 \qquad =\ (-1)^p\ (\psi g_p\smallsmile \varphi g_1) (f_{p+1} (1\otimes y^2 x^p\otimes 1)).
\end{array}
\]
The last case to consider is  $p,q\geq 2$, evaluating in $1\otimes x^{p+q}\otimes 1$ gives
\[
\begin{array}{rcl}
  (\psi \circ g_p \smallsmile \varphi \circ g_q) (f_{p+q} (1\otimes x^{p+q}\otimes 1)) &=& (\psi \circ g_p \smallsmile \varphi \circ g_q) (1\otimes x^{\otimes p+q}\otimes 1)\\
  &=& \varphi(g_q(1\otimes x^{\otimes q}\otimes 1)) \psi(g_p(1\otimes x^{\otimes p}\otimes 1)) \\
  &=& \varphi(1\otimes x^q\otimes 1) \psi(1\otimes x^p\otimes 1),
  \end{array}
  \]
  \[
  \begin{array}{l}
  (\psi \circ g_p \smallsmile \varphi \circ g_q) (c_{p,q} (f_{p+q} (1\otimes x^{p+q}\otimes 1))) = (\psi \circ g_p \smallsmile \varphi \circ g_q) (c_{p,q} (1\otimes x^{\otimes p+q}\otimes 1))\\
  \qquad \qquad \qquad \qquad \qquad \qquad \qquad =\ (\psi \circ g_p \smallsmile \varphi \circ g_q) (c_{p,q} (1\otimes x^{\otimes p}\otimes x^{\otimes q}\otimes 1))\\
  \qquad \qquad \qquad \qquad \qquad \qquad \qquad =\  (\psi \circ g_p \smallsmile \varphi \circ g_q) (1\otimes t^p\cdot (x^{\otimes q})\otimes x^{\otimes p}\otimes 1)\\
  \qquad \qquad \qquad \qquad \qquad \qquad \qquad =\ (\psi \circ g_p \smallsmile \varphi \circ g_q) (1\otimes (t^p\cdot x)^{\otimes q}\otimes x^{\otimes p}\otimes 1)\\
  \qquad \qquad \qquad \qquad \qquad \qquad \qquad =\ (\psi \circ g_p \smallsmile \varphi \circ g_q) (1\otimes ((-1)^p x)^{\otimes q}\otimes x^{\otimes p}\otimes 1) \\
  \qquad \qquad \qquad \qquad \qquad \qquad \qquad =\ (-1)^{pq}(\psi \circ g_p \smallsmile \varphi \circ g_q) (1\otimes x^{\otimes p+q}\otimes 1)\\
  \qquad \qquad \qquad \qquad \qquad \qquad \qquad =\ (-1)^{pq} (\psi \circ g_p \smallsmile \varphi \circ g_q) (f_{p+q} (1\otimes x^{p+q}\otimes 1)),
 \end{array}
 \]
 For $1\otimes y^2 x^{p+q-1}\otimes 1$ we get
\[
 \begin{array}{l}
 (\psi g_p\smallsmile \varphi g_q) (f_{q+p} (1\otimes y^2x^{q+p-1}\otimes 1)) \\
 =\ (\psi g_p\smallsmile \varphi g_q)(1\otimes y\otimes yx\otimes x^{\otimes q+p-2}\otimes 1 \\
 \qquad \quad -\ 1\otimes x\otimes (y^2+yx)\otimes x^{\otimes q+p-2}\otimes 1 \\
 \qquad \quad +\ \sum_{i=0}^{q+p-3}(-1)^i\ 1\otimes x^{\otimes 2+i}\otimes(y^2+yx)\otimes x^{\otimes q+p-3-i}\otimes 1)\\ 
    = \varphi g_q(1\otimes y\otimes yx\otimes x^{\otimes q-2}\otimes 1) \psi g_p(1\otimes x^{\otimes p}\otimes 1) \\ \quad -\ \varphi g_q(1\otimes x\otimes (y^2+yx)\otimes x^{\otimes q-2}\otimes 1) \psi g_p(1\otimes x^{\otimes p}\otimes 1)\\ 
    \quad +\ \sum_{i=0}^{q-3} (-1)^i \ \varphi g_q(1\otimes x^{\otimes 2+i}\otimes (y^2+yx)\otimes x^{\otimes q-3-i}\otimes 1) \psi g_p(1\otimes x^{\otimes p}\otimes 1)\\ \quad +\ \sum_{i=q-2}^{q+p-3} \varphi g_q(1\otimes x^{\otimes q}\otimes 1) \psi g_p (1\otimes x^{\otimes 2+i-q}\otimes (y^2+yx)\otimes x^{\otimes q+p-3-i}\otimes 1) \\
    = \varphi(1\otimes y^2x^{q-1}\otimes 1) \psi(1\otimes x^p\otimes 1) + (-1)^{q-2} \varphi(1\otimes x^q\otimes 1) \psi(1\otimes y^2x^{p-1}\otimes 1).
\end{array}
\]
Most of the terms in the above sums vanish because 
\[g_n(1\otimes x^{\otimes j}\otimes y^2\otimes x^{\otimes n-1-j}\otimes 1) = g_n(1\otimes x^{\otimes j}\otimes yx\otimes x^{\otimes n-1-j}\otimes 1) = 0, \ \forall n\geq2, j\geq 1.
\]
Regarding the unvanishing term, which corresponds to $i=q-2$, we have that 
\[
g_q(1\otimes y^2\otimes x^{\otimes q-1}\otimes 1) = 1\otimes y^2x^{q-1}\otimes 1, \hbox{ and }
g_q(1\otimes yx\otimes x^{\otimes q-1}\otimes 1) = y\otimes x^{q}\otimes 1,
\]
which vanishes after applying the $A$-bimodule morphism $\psi$. Finally
\[
\begin{array}{l}
   (\psi \circ g_p\smallsmile \varphi \circ g_q) (c_{p,q} (f_{q+p} (1\otimes y^2x^{q+p-1}\otimes 1))) \\
   =\ (\psi \circ g_p\smallsmile \varphi \circ g_q) (c_{p,q} (1\otimes y\otimes yx\otimes x^{\otimes p-2+q}\otimes 1 \\
   \qquad \quad -\ 1\otimes x\otimes (y^2+yx)\otimes x^{\otimes p-2+q} \otimes 1\\
   \qquad \quad +\ \sum_{i=0}^{p-3}(-1)^i\ 1\otimes x^{\otimes 2+i}\otimes(y^2+yx)\otimes x^{\otimes p-3-i+q}\otimes 1\\
   \qquad \quad +\ \sum_{i=p-2}^{q+p-3} (-1)^i\ 1\otimes x^{\otimes p}\otimes x^{\otimes 2+i-p}\otimes (y^2+yx)\otimes x^{\otimes q+p-3-i}\otimes 1))\\
    = (\psi \circ g_p\smallsmile \varphi \circ g_q) (1\otimes t^{p+1}(x^{\otimes q})\otimes y\otimes yx\otimes x^{\otimes p-2} \otimes 1 \\
    \qquad \quad -\ 1\otimes t^{p+1}(x^{\otimes q})\otimes x\otimes (y^2+yx)\otimes x^{\otimes p-2}\otimes 1\\
    \qquad \quad +\ \sum_{i=0}^{p-3}(-1)^i\ 1\otimes t^{p+1}(x^{\otimes q})\otimes x^{\otimes 2+i}\otimes(y^2+yx)\otimes x^{\otimes p-3-i}\otimes 1\\
    \qquad \quad +\ \sum_{i=p-2}^{q+p-3} (-1)^i\ 1\otimes t^p(x^{\otimes 2+i-p}\otimes (y^2+yx)\otimes x^{\otimes q+p-3-i})\otimes x^{\otimes p}\otimes 1) \\
        = (\psi \circ g_p\smallsmile \varphi \circ g_q) (1\otimes ((-1)^{p+1}x)^{\otimes q}\otimes y\otimes yx\otimes x^{\otimes p-2} \otimes 1 \\
    \qquad \quad -\ 1\otimes ((-1)^{p+1}x)^{\otimes q}\otimes x\otimes (y^2+yx)\otimes x^{\otimes p-2}\otimes 1\\
    \qquad \quad +\ \sum_{i=0}^{p-3}(-1)^i\ 1\otimes ((-1)^{p+1}x)^{\otimes q}\otimes x^{\otimes 2+i}\otimes(y^2+yx)\otimes x^{\otimes p-3-i}\otimes 1 \\
\qquad \quad +\ \sum_{i=p-2}^{q+p-3} (-1)^i\ 1\otimes ((-1)^{p}x)^{\otimes 2+i-p}\otimes (y^2- (p-1)yx-pxy) \\
\qquad \qquad \qquad \qquad \qquad  \otimes ((-1)^{p}x)^{\otimes q+p-3-i}\otimes x^{\otimes p}\otimes 1) 
\end{array}
\]
\[
\begin{array}{l}
= (-1)^{(p+1)q}\varphi (g_q(1\otimes x^{\otimes q}\otimes 1)) \psi (g_p (1\otimes y\otimes yx\otimes x^{\otimes p-2}\otimes 1)) \\
\quad -\ (-1)^{(p+1)q}\varphi (g_q(1\otimes x^{\otimes q}\otimes 1))  \psi (g_p (1\otimes x\otimes (y^2+yx)\otimes x^{\otimes p-2}\otimes 1))\\
\quad +\ (-1)^{(p+1)q} \sum_{i=0}^{p-3}(-1)^i \varphi (g_q (1\otimes x^{\otimes q}\otimes 1))  \\
\qquad \qquad \qquad \qquad \qquad \psi (g_p (1\otimes x^{\otimes 2+i}\otimes(y^2+yx)\otimes x^{\otimes p-3-i}\otimes 1)) \\
\quad +\ (-1)^{p(q-1)} \sum_{i=p-2}^{q+p-3}(-1)^{i} \varphi (g_q (1\otimes x^{\otimes 2+i-p}\otimes (y^2-(p-1)yx-pxy)\\
\qquad \qquad \qquad \qquad \qquad \qquad \qquad \qquad \otimes x^{\otimes q+p-3-i}\otimes 1))  \psi (g_p (1\otimes x^{\otimes p}\otimes 1))\\
= (-1)^{(p+1)q}\varphi(1\otimes x^q\otimes 1)\psi(1\otimes y^2x^{p-1}\otimes 1) \\
\quad +\ (-1)^{p(q-1)}(-1)^{p-2} \varphi(1\otimes y^2x^{q-1}\otimes 1) \psi(1\otimes x^p\otimes 1)\\
   = (-1)^{pq}\left((-1)^q\varphi(1\otimes x^q\otimes 1)\psi(1\otimes y^2x^{p-1}\otimes 1) + \varphi(1\otimes y^2x^{q-1}\otimes 1) \psi(1\otimes x^p\otimes 1)\right)\\
   = (-1)^{pq}(\psi \circ g_p\smallsmile \varphi \circ g_q) (f_{q+p} (1\otimes y^2x^{q+p-1}\otimes 1)).
\end{array}
\]
The second term, as well as the terms in the sums corresponding to $i\leq p-3$, vanish after applying $g_p$ in their second factor. Those corresponding to $i\geq p-2$ can be expressed as a sum of three terms by distributing over the expression in the tensor $y^2-(p-1)yx-pxy$. The only non vanishing one corresponds to $y^2$.

\bibliographystyle{alpha}
\bibliography{bibmono}

\end{document}